\documentclass[aip,pre]{revtex4-1}
\usepackage{graphicx,calc,empheq,amssymb}
\usepackage{amsmath,amsthm}
\usepackage{amsfonts}
\usepackage{hyperref}
\usepackage{bm}
\usepackage{mathrsfs}
\usepackage{amsmath}
\usepackage{amsthm}
\usepackage{eucal}  % Euler script symbols in math mode. Comando: \mathcal{A} (anziché \CMcal{A})
\usepackage{enumitem}  % enu­mer­ate, item­ize and de­scrip­tion
\usepackage[T1]{fontenc} %en­able hy­phen­ation to op­er­ate on texts con­tain­ing any char­ac­ter in the font
\usepackage[subnum]{cases}
\DeclareFontFamily{U}{mathc}{}
\DeclareFontShape{U}{mathc}{m}{it}%
{<->s*[1.03] mathc10}{}
\DeclareMathAlphabet{\mathscr}{U}{mathc}{m}{it}
\usepackage{color}
\usepackage{gensymb}  % definisce direttamente \degree e latri simboli, senza dover usare siunitx (vedi sotto)
\if@mathematic
   \def\vect#1{\ensuremath{\mathchoice
                     {\mbox{\boldmath$\displaystyle\mathbf{#1}$}}
                     {\mbox{\boldmath$\textstyle\mathbf{#1}$}}
                     {\mbox{\boldmath$\scriptstyle\mathbf{#1}$}}
                     {\mbox{\boldmath$\scriptscriptstyle\mathbf{#1}$}}}}
\else
   \def\vect#1{\ensuremath{\mathchoice
                     {\mbox{\boldmath$\displaystyle#1$}}
                     {\mbox{\boldmath$\textstyle#1$}}
                     {\mbox{\boldmath$\scriptstyle#1$}}
                     {\mbox{\boldmath$\scriptscriptstyle#1$}}}}
\fi

\if@vprodwedge
  
\else
  
\fi

{\bf}{\it}
{\bf}{\it}
{\bf}{\it}
{\bf}{\it}
{\bf}{\it}
{\bf}{\it}
\theoremstyle{lemma}
{\bf}{\bf}

\draft 
\begin{document}

        % Use the \preprint command to place your local institutional report number 
% on the title page in preprint mode.
% Multiple \preprint commands are allowed.
%\preprint{}
\title{ 
Colored Stochastic Multiplicative Processes with Additive Noise Unveil a Third-Order PDE, Defying Conventional FPE and Fick-Law Paradigms.
 }
\author{Marco Bianucci}
\email[]{marco.bianucci@cnr.it}
%\homepage[]{Your web page}
%\thanks{}
\affiliation{Istituto di Scienze Marine, Consiglio Nazionale delle Ricerche (ISMAR - CNR),\\
19032 Lerici (SP), Italy}
\author{Mauro Bologna}
%\email[]{marco.bianucci@cnr.it}
%\homepage[]{Your web page}
%\thanks{}
\affiliation{Departamento de Ingenier\'ia El\'ectrica - Electr\'onica, Universidad de Tarapac\'a,
Arica, 1000000, Chile}
\author{Riccardo Mannella}
\affiliation{Dipartimento di Fisica, Universit\`a di Pisa, 56100 Pisa, Italy}
\date{\today}
%       
        %Collaboration name if desired (requires use of superscriptaddress
        %option in \documentclass). \noaffiliation is required (may also be
        %used with the \author command).
        %\collaboration can be followed by \email, \homepage, \thanks as well.
        %\collaboration{}
        %\noaffiliation

\begin{abstract}
Research on stochastic differential equations (SDE) involving both additive and multiplicative noise has been extensive. 
In situations where the primary process is driven by a multiplicative stochastic process, additive white noise typically represents an 
intrinsic and unavoidable fast factor, including phenomena like thermal fluctuations, inherent uncertainties in measurement processes, or rapid
wind forcing in ocean dynamics. This work focuses on a significant class of such systems, particularly those characterized by linear drift and 
multiplicative noise, extensively explored in the literature.
Conventionally, multiplicative stochastic processes are also treated as white noise in existing studies. However, when considering colored 
multiplicative noise, the emphasis has been on characterizing the far tails of the probability density function (PDF), regardless of the 
spectral properties of the noise. In the absence of additive noise and with a general colored multiplicative SDE, standard perturbation approaches lead to a second-order PDE known as the Fokker-Planck Equation (FPE), 
consistent with Fick's law.
This investigation unveils a notable departure from this standard behavior when introducing additive white noise. At the leading order of the 
stochastic process strength, perturbation approaches yield a \textit{third-order PDE}, irrespective of the white noise intensity. The breakdown
of the FPE further signifies the breakdown of Fick's law. Additionally, we derive the explicit solution for the 
equilibrium PDF corresponding to this third-order PDE Master Equation. Through numerical simulations, we demonstrate significant deviations 
from outcomes derived using the FPE obtained through the application of Fick's law

%Research on stochastic differential equations (SDE) with additive and multiplicative noise has been extensive, particularly 
%in scenarios where the primary process is driven by a multiplicative stochastic process. 
%In this context, additive white noise represents an intrinsic and unavoidable fast factor. 
%This study focuses on systems with linear drift and multiplicative functions, a well-explored class in the literature.
%While multiplicative stochastic processes are typically treated as white noise, the emphasis shifts when considering colored 
%multiplicative noise. In the absence of additive noise, standard perturbation approaches yield a second-order PDE, i.e. a Fokker-Planck 
%Equation (FPE), consistent with Fick's law for general colored multiplicative SDEs.
%However, introducing additive white noise leads to a surprising departure from this standard behavior. At the leading order of the stochastic 
%process strength, perturbation approaches result in a \textit{third-order PDE}, regardless of white noise intensity. This breakdown of the FPE 
%also indicates the breakdown of Fick's law. We derive the explicit solution for the equilibrium PDF corresponding to this third-order PDE 
%Master Equation and demonstrate significant deviations from outcomes obtained using the FPE under Fick's law through numerical simulations.
\end{abstract}
\pacs{}% insert suggested PACS numbers in braces on next line

\maketitle %\maketitle must follow title, authors, abstract and \pacs
%+Contents
%-Contents
\section{Introduction}
\begin{figure}[htb]
\centering
\includegraphics[scale=.5]{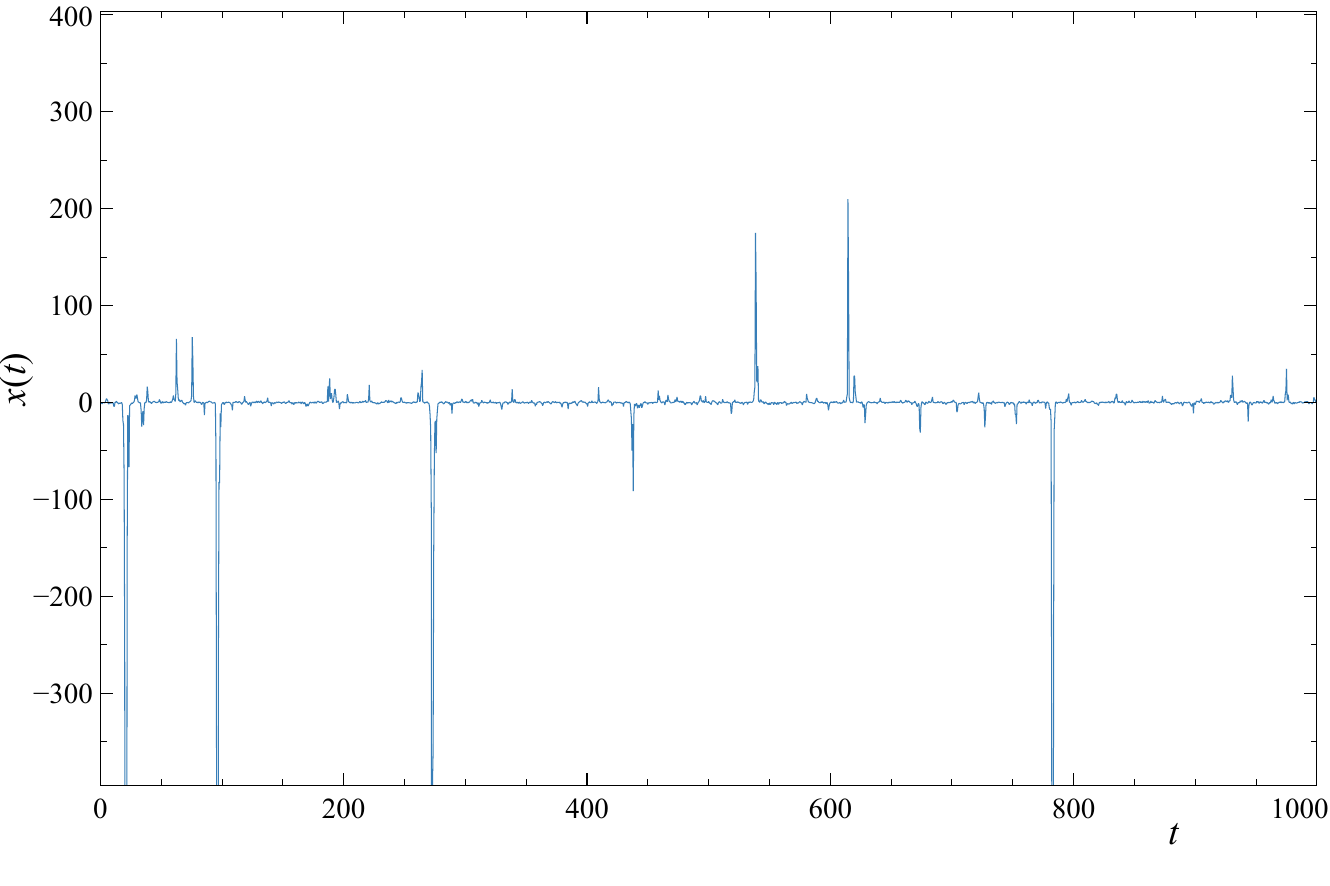}
\caption{A representative example, illustrating the intermittent behavior, 
is depicted in the time evolution of the amplitude \(x(t)\) for the SDE \eqref{SDE2_1} 
with parameters \(\tau=0.5\), \(\epsilon=5.0\), \(\gamma=2.0\), \(D_f =0.5\).
}
\label{fig:intermittenza}
\end{figure}

Linear equations forced by both additive and multiplicative noises are prevalent in
almost every scientific discipline. In the general $N$-dimensional case ($N$-D), these equations reads 
\begin{align}
\label{SDE2_MD_}
\dot {\vect x}= -\mathbb{E}\cdot \vect x+\vect f(t)
-    \bm{\Xi}(t)\cdot \vect x
\end{align}
where $\bm x:=(x_1,\dots, x_N$), $\mathbb{E}$ and $\bm{\Xi}(t)$ are $N\times N$ matrices with constant  %$\mathbb{E}_{i,j}%=e_{i,j}$ 
and stochastic components% 
%$\bm{\Xi}_{i,j}(t)=\xi_{i,j}(t)$
, respectively%
%and $i,j\in\{1,2,...,N\}$
; $\vect f(t)$ is a multidimensional white noise with correlation, or diffusion matrix given by $\mathbb{D}$.
As shown in\cite{dPA208b}, the extension to an infinite (or continuous) vector space of  \eqref{SDE2_MD_}, 
leads to a general model that describes a large class of important physical phenomena in fluid dynamics and  in quantum mechanics.
More in general, the model \eqref{SDE2_MD_} represents a random multiplicative process (RMP), a well-known mechanism
that gives rise to power-law behaviors.  
Widely employed as a model in various systems with both discrete and continuous time, the RMP
has been applied to phenomena such as on-off intermittency\cite{nPRE58,fyPTP74,yfPLA8,Pikovsky1992,phhPRL72}
and general intermittency (see Fig.~\ref{fig:intermittenza}) with power law statistics\cite{sPRE57,scJPII},
lasers\cite{sbPRA20,ghsPRL48}, economic activity\cite{lsIJMPC07,tstPRL79}, fluctuations in biological
populations within changing environments\cite{tTPB12}, and the advection of passive scalar fields by
fluids\cite{dPA208,dPA208b}. It  is clearly a paradigmatic model for theories
on large fluctuations (e.g., \cite{lCSF139} and references therein).

Therefore, the significance of the model \eqref{SDE2_MD_} cannot be overstated.

For simplicity, this work focuses on the 1-D version of the model  \eqref{SDE2_MD_}:
\begin{align}
\label{SDE2_1}
\dot {x}= -\gamma x+f(t)-  \epsilon\, x \xi (t),
\end{align}
which is the primary focus of most of the literature cited above. The extension to the $N$-D case \eqref{SDE2_MD_} of the formal results 
is straightforward yet somewhat intricate, as detailed in Appendix~\ref{app:multiD}. 

Thus, in \eqref{SDE2_1}  $f(t)$ is a white noise with diffusion coefficient $D_f$,  $\xi(t)$ is a Gaussian stochastic process with zero average,  finite correlation time $\bar \tau$\cite{Note1}
%\footnote{The general prescription is that there is a 
%        time $\bar\tau$ such that, for any
%        time $t$, the instances of $\xi$ at times $t'>t+\bar\tau$  are ``almost statistically uncorrelated'' with the instances of $\xi$ at  times $t'<t$. 
%        For ``almost statistically uncorrelated'' we mean that the joint probability density functions factorize up to terms $O(\bar\tau)$: $p_n(\xi ,t_1';\xi_2,t_2';...;\xi_k,t_k';\xi_{k+1},t_{1};...;\xi_n,t_h)=
%        p_k(\xi ,t_1';\xi_2,t_2';...;\xi_k,t_k')\,p_h(\xi_{k+1},t_{1};...;\xi_n,t_h)+O(\bar{\tau})$ with $k,h,n\in \mathbb{N}$,  $k+h=n$ and $t_i'> t_j+\bar{\tau}$.
%        For example, $p_2(\xi ,t';\xi_2,t)=p_1(\xi ,t')\,p_1(\xi_2,t)+O(\bar{\tau})$.\label{note1}}
%        
and normalized autocorrelation function  
$\varphi(t)=\langle \xi (t)\xi(0)\rangle_\xi/\langle\xi^2\rangle_\xi$. 
We  use the notation    $\langle ... \rangle_\xi$
to indicate the average over the realizations of the random process $\xi(t)$, which is assumed at equilibrium.  
We also define $\tau:=\int_0^\infty \varphi(u) \text{d}u$. The value of this integral can be much smaller than the decorrelation time $\bar \tau$,
in those cases when the function $\varphi(u)$ decays oscillating with time. 
Without  loss of generality, we assume that $\langle \xi^2\rangle_\xi=1$. 
Thus, the intensity of the fluctuations in the stochastic perturbation is governed by $\epsilon$. 
However, as demonstrated in Appendix~\ref{app:A}, the effective adimensional perturbation 
strength is measured by the parameter $\delta:=\epsilon \tau$. 
Strictly speaking, a more appropriate definition of $\delta$ should be $\delta=\epsilon \bar{\tau}$. 
This is because the real parameter that measures the perturbation strength in the 
approach illustrated in Appendix~\ref{app:A}, resulting in a series of cumulants, 
involves the correlation parameter $\bar{\tau}$ (as properly defined in note~\cite{Note1}) 
instead of $\tau$. 
Thus, for a more rigorous treatment, we should replace $\delta$ with 
$\delta \tau/\bar \tau$ in all the analytical results presented hereafter. 
However, to avoid complicating the formal expressions, we have chosen to stick 
with the current definition of $\delta$.

The drift field $-\gamma x$ in 
the SDE \eqref{SDE2_1}  can also be interpreted as originating from the same multiplicative stochastic process, when
its average is equal to $-\gamma/\epsilon$.
If $x(t) $ is intended
as the velocity of a Brownian particle, the SDE \eqref{SDE2_1} has  the important characteristic that it can be considered as
a continuous process realization of  L\'evy flights~\citep{lhuEPJB78} for some parameter range. 

We will refer to $f(t)$ in \eqref{SDE2_1} as the internal or intrinsic noise. 
This terminology is apt as it typically originates from intrinsic and unavoidable 
factors such as thermal fluctuations, inherent uncertainty in measurement processes, 
or rapid wind forcing in the context of ocean dynamics, among other possibilities.

As mentioned, the stochastic differential equation (SDE) \eqref{SDE2_1} has been extensively studied in the scientific literature. 
However, in almost all of these works, besides the additive noise $f(t)$, the stochastic process $\xi(t)$ has been considered as white noise. 
In cases where a colored stochastic process $\xi(t)$ has been considered, as in~\cite{dPA208b}, the focus of the work has been on 
characterizing the far tail of the probability density function (PDF) of $x$, which, as we will see hereafter, does not depend on 
the spectral properties of the multiplicative noise.

To the best of our knowledge, there are no papers that obtain a simple closed simple form (i.e., not a formal result with infinite 
series of operators) for the equation of the PDF in all its support range and the corresponding equilibrium solution. 
We will remedy this gap, focusing on the case in which the $\delta$ parameter is small, and we will find surprisingly 
simple results, but  not fitting either the Fokker-Planck equation (FPE) structure or Fick's law.

For an easy start, let us first assume that the intrinsic noise is absent.(i.e., $f(t)=0$). In this case it is easy to show that, {\em regardless of the values of $\tau$ and $\epsilon$}, 
 the Master Equation (ME) for the PDF of $x$ in \eqref{SDE2_1} coincides with the following 
 FPE (see Appendix~\ref{app:A} for simplicity, we use the shorthand 
$\partial_y:= \partial/\partial y$):
\begin{align}  
\label{FPE_Linear} 
&\partial_t {P}(x;t)=
\left\{\gamma \partial_x x
+ \frac{\delta^2}{\tau}\partial_x x
\partial_x\,x
\right\} P(x;t).
\end{align}  
This fact indicates that the process 
\eqref{SDE2_1}, with $D_f=0$, does not depend on the spectrum (or color) features of the stochastic process $\xi(t)$. 
Consistently, in the white noise limit, i.e., for $\tau\to 0$ and $\delta^2/\tau=\epsilon^2\tau$ held constant, the FPE~\eqref{FPE_Linear} remains unaltered and 
 corresponds to the standard FPE for SDEs with multiplicative white noise, under the 
Stratonovich interpretation of Wiener process differentials. 

The FPE~\eqref{FPE_Linear} can also be rewritten as a conservative equation as
\begin{align}  
\label{FPE_Linear_Fick} 
\partial_t {P}(x;t)=\partial_x J(x)
\end{align}  
with
\begin{align}  
\label{J}
J(x):=\left\{(\gamma\tau +\delta^2) x/\tau+ D_\xi(x)
\partial_x
\right\} P(x;t).
\end{align}  
where we have introduced the inhomogeneous diffusion coefficient, 
\begin{equation}
\label{Dxi}
    D_\xi(x):= \delta\epsilon x^2=\delta^2x^2/\tau,
\end{equation} 
characteristic of multiplicative noise.
The interpretation of the terms appearing in Eq.~\eqref{J} is straightforward: 
in absence of internal noise ($D_f=0$ in \eqref{SDE2_1}), the  
multiplicative stochastic process generates an additional 
friction/drift term proportional to the intensity of the stochastic perturbation and an inhomogeneous diffusion process, 
proportional to the gradient of the PDF (thus, following Fick's law). 
The  equilibrium PDF of 
\eqref{FPE_Linear_Fick} is obtained by setting $J(x)=0$ in \eqref{J}, yielding $P(x)_{eq}\propto x^{-\left(1+\frac{\gamma\tau}{\delta^2}\right)}$,  showcasing  a singular behaviour (it is non-integrable)
 around $x=0$\footnote{The cause of this issue lies in the fact that $x=0$ is both a stable attraction
point for the unpertubed motion 
($\epsilon \tau^2=0$)  and
a point where the perturbation vanishes. Hence, given that any perturbed
trajectory will eventually arrive at 
$x=0$, this is an accumulation point of any initial ``ensemble''.}. 

The introduction of an internal diffusion source effectively addresses this issue and is physically plausible for many realistic models. 
In fact, by leveraging Fick's law and including the standard, constant diffusion coefficient $D_f$ in the current \eqref{J},
we have
\begin{align}  
\label{J:Fick}
J(x):=\left\{(\gamma\tau +\delta^2) x/\tau+ (D_\xi(x)
+ D_f )\partial_x
\right\} P(x;t).
\end{align}  
The equilibrium PDF, obtained by setting $J(x)=0$, is 
$P(x)_{eq}\propto 
(D_f+D_\xi(x))^{-\frac{1}{2}\left(1+\frac{\gamma\tau}{\delta^2}\right)}$. 
which no longer displays a singular behavior around $x=0$. 

Note that for $D_\xi(x)\gg D_f$, i.e., for $x\gg \sqrt{D_f\tau/\delta^2}$, it behaves similarly to the previous noiseless case. 
Thus, if there are no boundary conditions that constrain the variable $x$ to a finite range, the existence condition for the moments of $x$ remains unaltered by the introduction of this diffusive term.

We can say that the white noise $f(t)$, corresponding to a diffusion process with a diffusion coefficient $D_f$, 
introduces a repulsion from the origin, preventing any path from getting trapped at the 
$x=0$ point once reached. 
While the introduction of such intrinsic noise in multiplicative processes has been 
acknowledged by many researchers (see the same works already
cited above), however it has been not highlighted the fact that
even though the two fluctuating processes  are assumed independent
of each other, in effect their contributions  to the current $J(x)$ of \eqref{FPE_Linear_Fick} 
{\em don't simply add up}, unless the multiplicative process is a white noise too.
More precisely, in this work we will show that we have
\begin{equation}
\label{J:noFick} 
J(x)=\left\{(\gamma\tau +\delta^2) x/\tau+ (D_\xi(x)+D_f)\partial_x 
+ D_f D_\xi(x) \vartheta \partial^2_x \right\}P(x;t)
\end{equation}
with $\vartheta$, given in \eqref{vartheta}, having the dimension of time and coinciding with 
$2\tau$ for $\gamma \tau\ll 1$ and with $\gamma^{-1}$ for 
$\gamma \tau\gg 1$. 
Thus, in this case {\em  the Fick's law and the corresponding
FPE structure break down}.

The last term on the right-hand side of \eqref{J:noFick}, which invalidates Fick's law, arises from two factors:  
the finite time scale of the external stochastic perturbation (colored noise) and the non-commutativity 
of its Liouvillian with the Liouvillian associated with the internal noise (the standard diffusion operator). 
Specifically, when averaging over the external stochastic process, 
the expansion in cumulants of the PDF coincides with a power series 
of the adimensional parameter $\delta$ (see Appendix~\ref{app:A}). 
In the second order (which is the leading one for weak perturbations), these two factors yield, 
in the ME of the PDF of $x$, a correction to the FPE obtained by applying Fick's law. 
This correction is proportional to both $\delta$ and $D_f$ and turns out to be a {\em third-order partial differential operator on $x$}.

This  result is confirmed by the numerical simulations reported in Section \ref{sec:OU}.
We will undertake a detailed derivation of this phenomenon in the next session.
However, it's crucial to highlight that this 
departure from the standard FPE/Fick's law is a general observation, applicable beyond the linear drift case of~\eqref{SDE2_MD_}, and always occurs when both additive noise and multiplicative colored stochastic processes are present. For simplicity, here we focus on the linear 1D case of \eqref{SDE2_1}, while a more in-depth exploration of these findings will be presented in future works.

\section{A third order PDE for the PDF\label{sec:third_order}}
Given the infinitely short time correlation of the additive noise $f(t)$, to any realization
$\xi(\cdot)$ of the stochastic process $\xi(u)$, $0\le u \le t$, from the 
SDE \eqref{SDE2_1} we can write 
  the following  Liouville equation  for the PDF of $x$: $ {P}_{\xi(\cdot)}(x,t) $:
\begin{align}
\label{Liouville}
&\partial_t {P}_{\xi(\cdot)}(x,t) \nonumber \\
&= 
\{{\cal L}_a+ \epsilon\, \xi(t) {\cal  L}_{I}\}  {P}_{\xi(\cdot)}(x,t) ,
\end{align}
in which  ${\cal L}_a$ is the  unperturbed Liouville operator given by 
\begin{equation}
\label{La}
{\cal L}_a:=\gamma \partial_x x+D_f\partial_x^2;
\end{equation}
%
%${\cal L}_b$ is the Liouvillian operator of the fluctuating perturbation $\xi(t)$ (that in principle
%could also be a chaotic deterministic dynamical variable), 
and $ \epsilon\, \xi(t) {\cal  L}_{I}$
is the Liouville perturbation operator with:
\begin{equation}
\label{LI}
{\cal L}_I:=\partial_x x.
\end{equation}
If the perturbing process $\epsilon \xi(t)$ is weak (characterized by small values of the $\delta$
parameter), applying a perturbation
projection/cumulant approach 
\citep{bgJCP96,bmJPC4,bbgwgPLA174,bmxgwPRE47,bwgPLA190,bPRE70,Zwanzig2001,grigolini1989,bmwgPRE51,bJSTAT2015,bbmSUB} to Eq.~\eqref{Liouville},
at the leading order of $\delta$, we formally obtain the following standard result for the reduced
PDF of $x$ (see see Appendix~\ref{app:A} for details, note that hereafter $P(x;t):=\langle{P}_{\xi(\cdot)}(x,t)\rangle_\xi$, where $\langle...\rangle_\xi$ is the
average over the realizations of the stochastic process $\xi(t)$; throughout this work we will consider  $t\gg\tau$):
\begin{equation}  
\label{FP_F2} 
\partial_t {P}(x;t)= {\cal L}_a P(x;t)+ \frac{\delta^2}{\tau^2}  {\cal L}_I\int_0^{\infty}\text{d}u\,   \varphi(u)\,\tilde{\cal L}_I(-u)\;{P} (x;t),
\end{equation}  
where $\varphi(t)$ is the normalized autocorrelation function
of $\xi(t)$, as defined in the Introduction, and
\begin{equation}
\label{LIt}
\tilde{\cal L}_I(t) := e^{-{\cal L}_a t}{\cal L}_I e^{{\cal L}_a t}
\end{equation}
 is the interaction representation of the perturbing Liouvillian ${\cal L}_I$. By exploiting
the Hadamard's lemma for exponentials of operators we can also write 
\begin{equation}
\label{LIt2}
\tilde{\cal L}_I(t) =e^{-{\cal L}_a^\times t} [{\cal L}_I]
\end{equation}
 in which, for any couple of operators ${\cal A}$ and ${\cal B}$, we have defined 
${\cal A}^\times[{\cal B}]:=[{\cal A},{\cal B}]={\cal A} {\cal B}- {\cal B}{\cal A}$.
In literature (e.g.~\cite{bJMP59}),   $e^{{\cal A}^\times t}\left[{\cal B}\right]$ is called 
the Lie evolution of the operator ${\cal B}$ 
along 
 ${\cal A}$, for a time $t$.  

Because the perturbing Liouvillian ${\cal L}_I$ of \eqref{LI} is a first order differential operator,
the order of the differential operator corresponding to the second addend in the r.h.s. of 
\eqref{FP_F2} is obtained by adding to one  the order of differential operator of 
$\tilde{\cal L}_I(-u)$.
From Eq.~\eqref{LIt2}, we see that this latter is  the result of the
Lie evolution of   ${\cal L}_I$
along the unperturbed Liouvillian ${\cal L}_a$.  

If the decay time of $\varphi(u)$ is significantly shorter than $1/\gamma$, 
we can safely assume the approximation 
$\tilde{\cal L}_I(-u)\approx {\cal L}_I$ inside
the integral on the r.h.s. of \eqref{FP_F2}. 
Consequently, the ME \eqref{FP_F2} effectively reduces to a FPE.
However, when this is not the case, we must address the challenge of evaluating the full Lie 
evolution of   ${\cal L}_I$
along  ${\cal L}_a$. 
In-depth exploration of this topic, from a formal and general perspective, can be found 
in \cite{bJMP59}.
Specifically, Proposition 1 in \cite{bJMP59} is of particular relevance. 
For the simple 1-D case with linear drift, corresponding to the present
SDE~\eqref{SDE2_1}, we can derive the Lie evolution of ${\cal L}_I$ 
along the unperturbed Liouvillian ${\cal L}_a$ in \eqref{LIt2} as follows. From \eqref{LIt2} we 
have
\begin{equation}
\label{dtLIt}
\frac{d}{dt}\tilde{\cal L}_I(t) =-{\cal L}_a^\times \left[e^{-{\cal L}_a^\times t} 
[{\cal L}_I]\right]=-e^{-{\cal L}_a^\times t} \left[[{\cal L}_a,{\cal L}_I]\right].
\end{equation}
By using \eqref{La} and \eqref{LI}, we get
\begin{equation}
\label{LalI}
\left[{\cal L}_a,{\cal L}_I\right]
=[\gamma \partial_x x+D_f\partial_x^2\,,\,\partial_xx]=2D_f\partial_x^2
\end{equation}
thus, Eq.~\eqref{dtLIt} can be written as
\begin{align}
\label{dtLIt2}
\frac{d}{dt}\tilde{\cal L}_I(t) =&-2e^{-{\cal L}_a^\times t} 
\left[D_f\partial_x^2\right]\nonumber \\
=&-2e^{-{\cal L}_a^\times t} 
\left[\gamma \partial_x x+D_f\partial_x^2-\gamma \partial_x x\right]
=-2(\gamma \partial_x x+D_f\partial_x^2)+2\gamma \tilde{\cal L}_I(t)
\end{align}
of which the solution is
\begin{equation}
\label{expLalinearRed}
\tilde {\cal L}_I(t)=
\partial_x\,x+ D_f \frac{1-e^{2 \gamma  t}}{\gamma }
\partial_x^2.
\end{equation}
By using Eq.~\eqref{expLalinearRed} into the ME \eqref{FP_F2}, and  exploiting again  \eqref{La} and \eqref{LI}, we finally obtain 
\begin{align}  
\label{FP_Final} 
\partial_t {P}(x;t)&=
\left\{\gamma \partial_x\, x +D_f\partial _x^2 
+ \frac{\delta^2}{\tau} \partial_x x\partial_x x +
D_f \delta^2 \frac{\vartheta}{\tau}\partial_x\,x
\partial_x^2
\right\} P(x;t)\nonumber \\
&=\partial_x J(x)
\end{align}  
with $J(x)$ is given in \eqref{J:noFick} and the time $\vartheta$ defined as
\begin{align}  
\label{vartheta} 
\vartheta:=\frac{1}{\gamma\tau}\left(
 \tau-\hat \varphi(2\gamma)\right).
\end{align}  
The hat over a function indicates its Laplace transform: 
$\hat \varphi(s) := \int_0^{\infty}\text{d}u\,\varphi(u) e^{-s  u}$.
  The third order PDE \eqref{FP_Final} with \eqref{vartheta} is the main result of the present work. At the leading order
in powers of the $\delta$ parameter, Eq.~\eqref{FP_Final} is exact, 
irrespective of the value of the
diffusion coefficient $D_f$. 

Thus, upon 
 introducing the internal noise, alongside 
 the standard diffusion process, \textit{an additional mutual contribution is activated}. 
 As we can observe from \eqref{FP_Final}, this mutual contribution of the white internal noise
and the external multiplicative stochastic process takes on an odd nature 
 in terms of partial derivatives. 
As previously emphasized in the Introduction, we reiterate that the time parameter $\vartheta$ of \eqref{vartheta} coincides with
$2\tau$ for $\gamma \tau\ll 1$ and with $\gamma^{-1}$ for
$\gamma \tau\gg 1$. Consequently, the adimensional parameter 
$r:=\epsilon \vartheta$ is akin the $\delta$ parameter, but is rescaled based on the time scale relationship between the stochastic process and the unperturbed dynamics.

Imposing the equilibrium condition to   the ME \eqref{FP_Final}, i.e., setting  
$J(x)=0$, we 
obtain two different  analytical solutions, both involving the   Kummer confluent 
hypergeometric function of first 
kind\footnote{$\,_1F_1\left(a;b;z\right):=\sum _{k=0}^{\infty }  \frac{a_k}{b_k k!}$, where 
$(x)_n=x (x+1) \ldots  (x+n-1)=\Gamma  (x+n)/\Gamma(x) $ is the Pochhammer symbol.}:
\begin{align}
\label{HYP_r}
P_1(x)=& \,  _1F_1\left(\frac{1}{2} \left(\frac{\gamma\tau }{\delta^2  }+1\right);\frac{1}{2} \left(\frac{1}{\delta r  }+1\right);-\frac{x^2}{2 D_f  \vartheta}\right)\nonumber \\
P_2(x)=&\, 2^{\frac{1}{2} \left(\frac{1}{\delta r}-1\right)} 
r^{\frac{1}{2} \left(\frac{1}{\delta r}-1\right)} D_f^{\frac{1}{2} 
\left(
\frac{1}{\delta r}-1
\right)} 
x^{1-\frac{1}{\delta r}}\nonumber \\
&\times \,_1F_1\left(
\frac{r \left(2 \delta\epsilon+\gamma \right)-1}{2 \delta r}
;\frac{3}{2}-\frac{1}{2 \delta r};-\frac{x^2}{2D_f  \vartheta}
\right).
\end{align}
This fact is due to the third order nature of the PDE \eqref{FP_Final}. From a mathematical
point of view any linear combination of these two functions is also a possible solution. 
However, it is easy to show that the second one is not physically acceptable.
In fact, let us consider the behaviour of these two functions around $x=0$. We have
\begin{align}
P_1(x)\approx&\, 1 -\frac{x^2 \left[\gamma\tau +  \delta^2\right]}
{(2D_f\tau) \left(\delta r+1\right)}+O\left(x^3\right)
\label{HYP_r_0}\\
P_2(x)\approx&\, 2^{\frac{1}{2} \left(\frac{1}{\delta r}-1\right)} r^{\frac{1}{2} 
\left(
\frac{1}{\delta r}-1
\right)} 
D_f^{\frac{1}{2} \left(\frac{1}{\delta r}-1\right)} 
x^{1-\frac{1}{\delta r}}+O\left(x^3\right).
\end{align}
We see that if $R:=\delta r<1$, a condition which is typically met, the solution $P_2(x)$ is not integrable, 
therefore it must be discarded. 
The expression of the function \(P_1(x)\) in \eqref{HYP_r_0} implies that the presence of $r>0$ 
smears the equilibrium PDF around $x=0$.

Thus, the final result is given by
\begin{align}
\label{HYP_eq}
P_{eq}(x)\propto \, _1F_1\left(\frac{1}{2}
\left(\frac{\gamma \tau}{\delta^2}+1\right)
;\frac{1}{2} \left(\frac{1}{\delta r }+1\right);-\frac{x^2}{2 D_f  \vartheta}\right).
\end{align}
In the case in which the support of the PDF is not limited (for example, if there are not reflecting boundary conditions at some finite values of $x$),
we can evaluate the asymptotic behavior as  $x\to \pm\infty$ of the equilibrium PDF 
in \eqref{HYP_eq}, and the result is
$P_{eq}(x)\sim \left| x\right| ^{-\left(\frac{\gamma \tau}{\delta^2}+1\right)}$.
From this expression we observe that even when considering the contribution from the third partial derivative, 
the far tails of the equilibrium PDF of $x$ remain unaffected by the presence of the additive white noise $f(t)$. 
This implies that, for the case of infinite support of the PDF, the condition for the existence of the 
$n$-th moment of $x$ depends only on the fraction $\gamma\tau/(\delta^2)$, i.e., it remains independent 
of the spectral properties of $\xi(t)$.
Still in the case of an unbounded domain for PDF, from the PDE \eqref{FP_Final} it is possible
to obtain the following time differential equation for the 
moment of order $n$ of $x$:
\begin{align}  
\label{moments_linear} 
\partial_t \langle x^n \rangle=&-n\gamma
\langle x^n \rangle\left(
1- n\delta^2/(\gamma\tau)\right)\nonumber \\
&+n(n-1)
D_f\left(1-n\delta  r\right)
  \langle x^{n-2} \rangle.
\end{align}
For any fixed $n$, Eq.~\eqref{moments_linear} is a linear relationship between the first $n$
moments of $x$. The eigenvalues of  the 
corresponding matrix are $-n\gamma\left(
1- n\delta^2/(\gamma\tau)\right)$, thus, 
they do not depend on $D_f$ and $r$. Therefore, the relaxation behaviour
of the moments is independent on $D_f$ and $r$ as well and they exist
only if $\left(1- n\delta^2/(\gamma\tau)\right)>0$.
On the other hand, it is clear from the same
Eq.~\eqref{moments_linear} that the equilibrium values of the moments (if they exist) do depend on $D_f$, and also
on the value of  
$R:=\delta r$. 
When it exists, the equilibrium solution of Eq.~\eqref{moments_linear} is given by 
\begin{equation}
\label{moments_linear_eq}
\langle x^n \rangle_{eq}=
\left\{
\begin{array}{lr}
0&\text{for $n$ odd}\\
\left(\frac{D_f }{\gamma }\right)^{n/2} (n-1)!!\prod_{j=1}^{n/2} 
\frac{
 \left(1-2 j \, \delta r\right)}
 {\left(1-2j \,\delta^2/(\gamma\tau)\right)}&\text{for $n$ even}
 \end{array}
  \right.
\end{equation}
From \eqref{vartheta} we always
have $\vartheta<1/\gamma$, from which 
$\left(1- 2j\delta^2/(\gamma\tau)\right)>\left(1-2 j \, \delta r\right)$.
Therefore, if 
$\left(1- n\delta^2/(\gamma\tau)\right)>0$ (as it must be for the $n$-th moment do exist),
then  un-physical situation of an even moment smaller than
zero is incompatible
with Eq.~\eqref{moments_linear_eq}. 

This observation appears to contradict both the earlier findings that around 
$x=0$ the equilibrium PDF  broadens with increasing 
$R$ and that the far tails of this PDF do not depend on $R$. The explanation of this apparent
contradiction lies in the fact as we move away from the 
origin, where the expansion \eqref{HYP_r_0} holds, but before we arrive at the
asymptotic tail, the equilibrium PDF decays more quickly as function of $x$ due to the presence of $R>0$. This fact is easily confirmed by plotting the equilibrium PDF with $R=0$
and $R\ne0$ (see next section).  
 
If there
is a very large time scale separation, i.e. $\gamma\tau\ll 1$, from \eqref{vartheta} we have 
\begin{align}  
\label{r_1} 
R=\delta r\approx 
\epsilon^2\int_0^{\infty}\text{d}u\,\varphi(u)u\approx 
2 (\epsilon \tau)^2=2\delta^2, 
\end{align}  
that does not depend on
$\gamma$. The reader should note that in the white noise limit, i.e., as
$\tau$ approaches zero while keeping $\epsilon^2\tau=\delta \epsilon$ fixed, $R$ in
\eqref{r_1} tends to zero. Consequently, the non-Fick contribution to $J$ becomes
negligible as well. Conversely, if $\delta$, the relevant parameter for the cumulant series, is held fixed (small enough to allow truncating the series to the second cumulant),
while changing the time scale of the noise, $R$ in \eqref{r_1} remains constant.

In essence, with a small strength of the stochastic process $\xi(t)$ 
(while it diverges in the white noise limit), the breakdown of Fick's law and of the associated FPE for the model \eqref{SDE2_1} persists, resulting in a ME with a third derivative term.

\section{The case of Ornstein Uhlenbeck external stochastic process: analytical 
and numerical results\label{sec:OU}}

To explore the full range of $\tau$ and $\gamma$, we consider the specific
case
of exponentially decaying  correlation function: $\varphi(u)=\exp(-u/\tau)$, from which, by 
also exploiting
\eqref{vartheta}, we have:
\begin{align}  
\label{r_linear} 
\vartheta=\frac{2  \tau }{(2 \gamma  \tau +1)} 
\end{align}  
i.e., 
\begin{align}  
\label{R_linear} 
R:=\delta r=
\frac{2\delta^2 }{(2 \gamma  \tau +1)}. 
\end{align}  
We observe that $R$  depends solely on $\delta$ (the small parameter in
the cumulant expansion) and $\gamma\tau$ (quantifying the time scale 
separation between the unperturbed relaxation process and the relaxation of the 
correlation function of $\xi(t)$). It is evident from \eqref{R_linear} that 
$R$ decreases when the time scale separation decreases ($\gamma \tau$ increases) and 
increases 
quadratically with $\delta$.
 %For example, let us assume
 %  $\gamma  \tau=\epsilon\tau =0.5$, then we get $R=1/4$ see figure~\ref{fig:D3_OU}. 
 %
 \begin{figure}[htb]
\flushleft
 \includegraphics[width=1.0\textwidth]{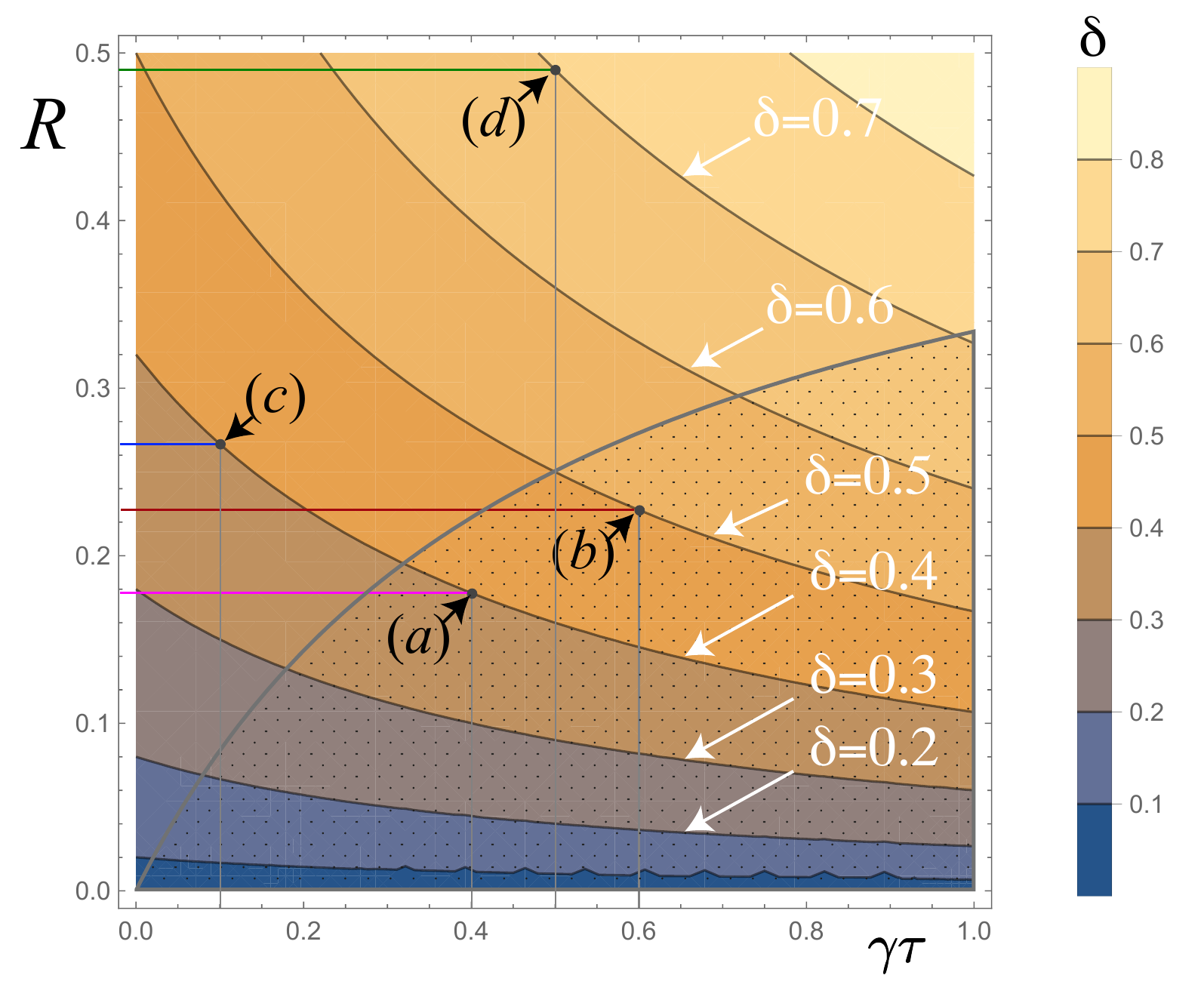}
 \caption{Plot of $R$ of \eqref{R_linear} vs $\gamma\tau$, for  various values of 
 $\delta$ (distinct curves). We see that at fixed $\delta$, as
$\gamma\tau$ decreases, $R$ increases. The same happens increasing  $\delta$, at
$\gamma\tau$ fixed.  The area with dotted background correspond to   $\delta$ and
$\gamma\tau$ values for which the variance of $x$ is finite ($\gamma\tau- 2\delta^2$>0). The
points in the graph labeled with the letters (a) and (b) ((c) and (d)) corresponds to the 
$\gamma\tau$ and $\delta$
values, used for the four plots of the PDF of figure~\ref{fig:PDF_OU1} (figure 
\ref{fig:PDF_OU2}).   }
 \label{fig:D3_OU}
\end{figure} 
The equilibrium PDF \eqref{HYP_eq} in this case reads:
\begin{align}
\label{HYP_OU}
P_{eq}(x)=N  \,_1F_1\left(\frac{1}{2} 
\left(\frac{\gamma \tau}{\delta^2}+1\right);\frac{1}{4} 
\left(\frac{2\gamma \tau +1}{\delta^2}+2\right);
-\frac{(2\gamma \tau +1) x^2}{4 D_f \tau }\right).
\end{align}
where $N$ is a normalization factor.
We note that, except for the quantity $D_f\tau$, which acts as a scale factor for $x$, also 
the equilibrium PDF~\eqref{HYP_OU} depends only on $\delta$ and $\gamma\tau$. 

\begin{figure}[htb]
\flushleft
%  \noindent
\includegraphics[width=1.0\textwidth]{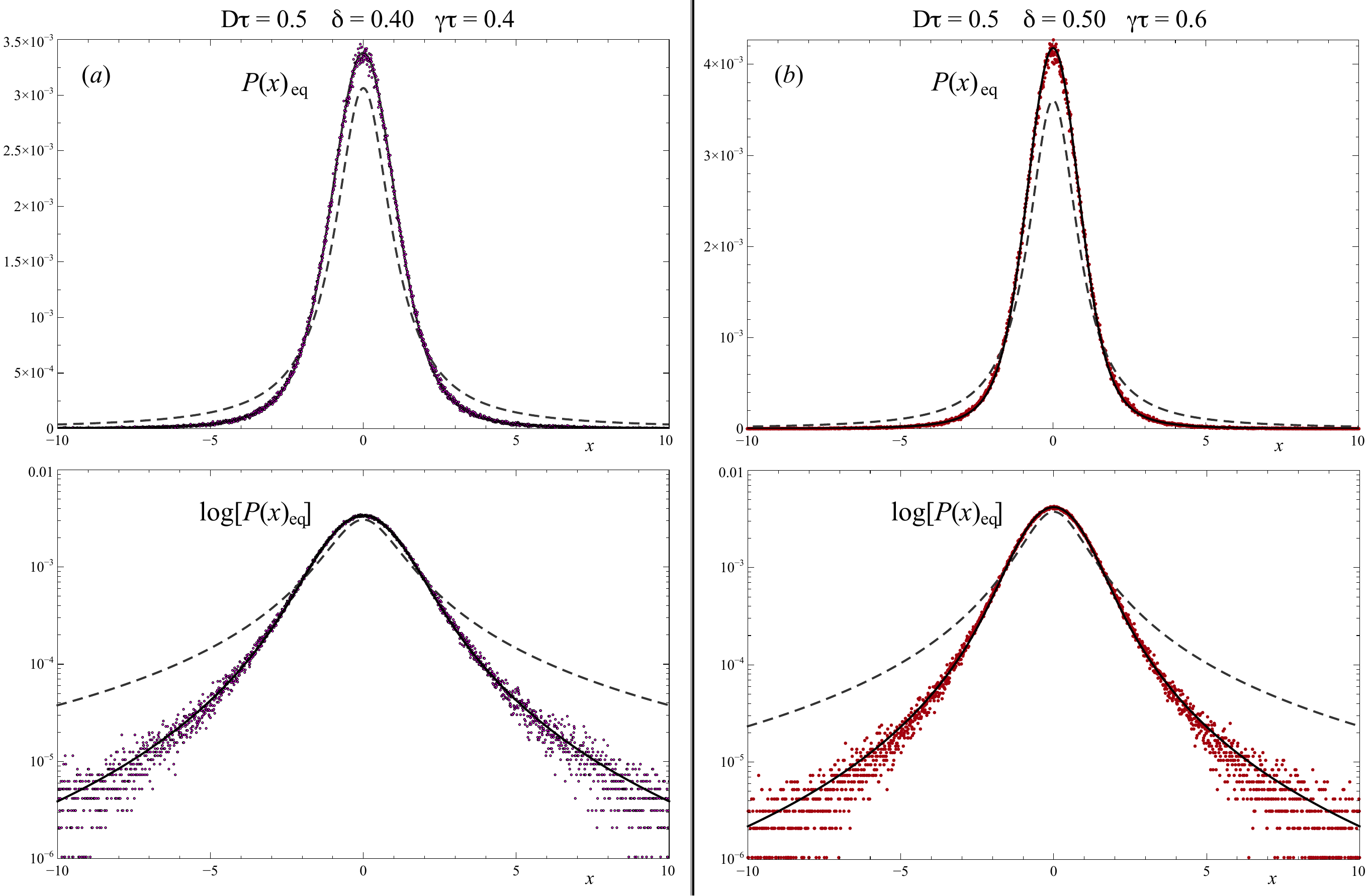}
\caption{Two vertical panels (a) and (b), respectively, displaying the 
Equilibrium PDF of the SDE~\eqref{SDE2_1} for the case in which the stochastic process
$\xi(t)$ is the Ornstein Uhlenbeck process. Panel (a): $D_f\tau=0.5$, $\delta=0.4$ and 
$\gamma\tau=0.4$.  Panel (b): $D_f\tau=0.5$, $\delta=0.5$ and $\gamma\tau=0.6$. 
In the bottom part, the same data as the upper part are presented in semi-logarithmic scale.  
Circles are the results of the numerical
simulation. 
Solid lines represent the theoretical result \eqref{HYP_OU}, i.e., the equilibrium solution
of the PDE \eqref{FP_Final}, with $r$ given in \eqref{r_linear}. Dashed lines depict
$P_{eq,FPE}(x)\propto (D_f+D_\xi(x))^{-\frac{1}{2}\left(1+\frac{\gamma\tau}{\delta^2}\right)}$,
the solution for the vanishing ``Fick''  current of \eqref{J:Fick} (or the equilibrium PDF of the corresponding FPE).
In these two cases,
corresponding to the two points (a) and (b) in the diagram of figure~\ref{fig:D3_OU},
we have $\gamma\tau>2\delta^2$, thus the variance of $x$ is finite (see text for details).}
\label{fig:PDF_OU1}
\end{figure}

\begin{figure}[htb]
\flushleft
5 \noindent
\includegraphics[width=\textwidth]{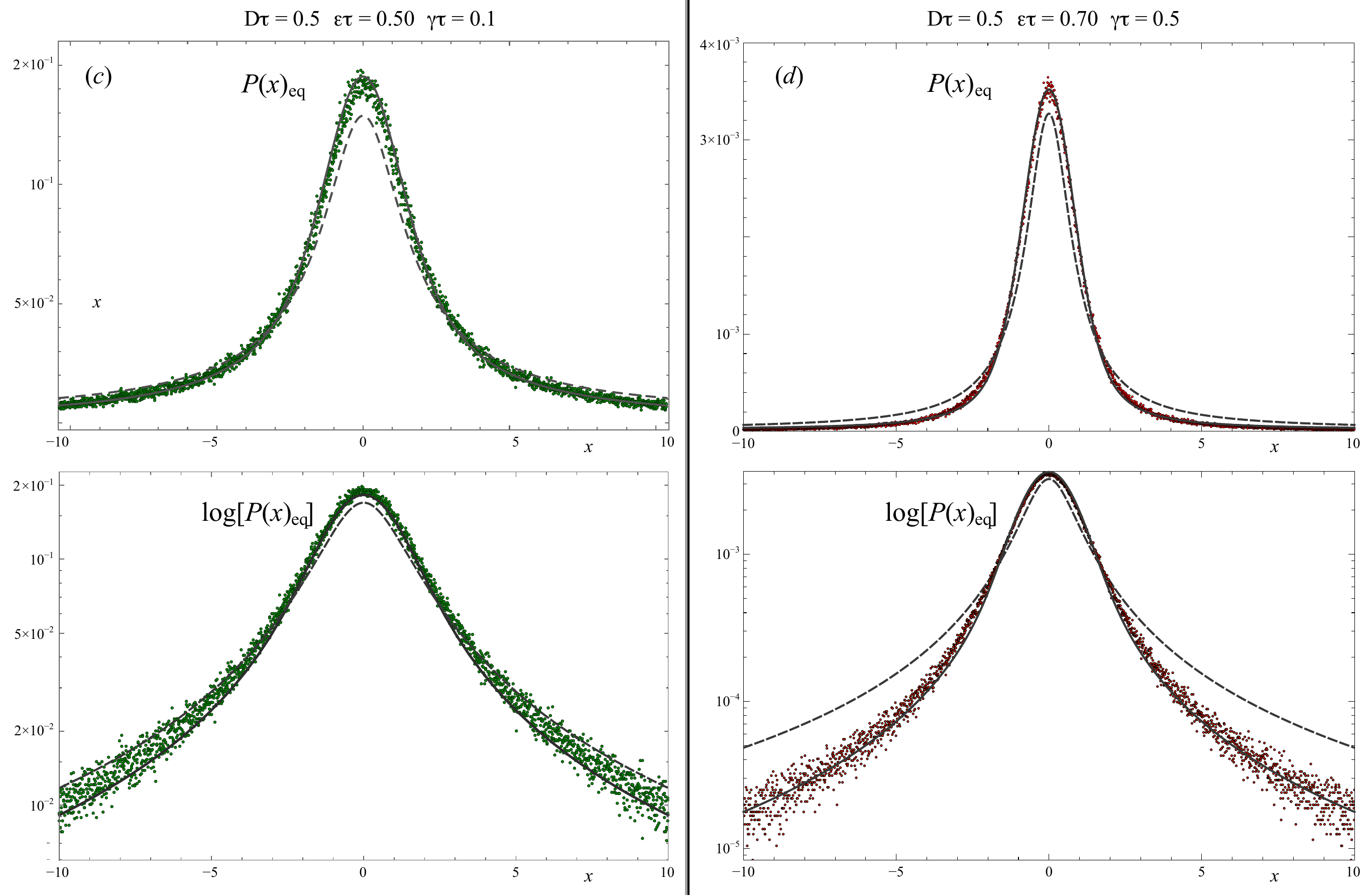}
\caption{The same as figure~\ref{fig:PDF_OU1}, but for different values  of $\gamma\tau$
and $\delta$ as indicated in the header of the panels. In these two cases,
corresponding to the two points (c) and (d) in the diagram of figure~\ref{fig:D3_OU},
we have $\gamma\tau<2\delta^2$, thus, at equilibrium all the moments of $x$ diverge.} 
\label{fig:PDF_OU2}
\end{figure}

In figures \ref{fig:PDF_OU1}-\ref{fig:PDF_OU2}, solid lines depict the plots of the
PDF~\eqref{HYP_OU} for a fixed $D_f\tau=0.5$ and different values of $\delta$ 
and $\gamma\tau$ corresponding to the points $(a)-(d)$ in the diagram \eqref{fig:D3_OU}. 
We have also included the corresponding results of the numerical
simulation of the SDE~\eqref{SDE2_1} (circles), where $\xi(t)$ is the 
Ornstein-Uhlenbeck process. 
Additionally, to assess the relevance of the non Fick 
contribution to the current, 
we have also plotted, with dashed lines, the normalized function
$P_{eq,FPE}(x)\propto
(D_f+D_\xi(x))^{-\frac{1}{2}\left(1+\frac{\gamma\tau}{\delta^2}\right)}$
which is the  solution for the vanishing ``Fick'' current of 
\eqref{J:Fick} (or the equilibrium PDF of the corresponding FPE).
The excellent agreement of the analytical result \eqref{HYP_OU} with numerical simulations 
is evident, while when relying on the $P_{eq,FPE}(x)$, 
the comparison with numerical simulations is not at all so good.
% For the numerical simulations and the comparison with the theoretical results, we consider
% the case where moov   $\xi(t)$
% is here a Ornstein Uhlenbeck (OU) external stochastic process. Thus, the $r$ parameter is
% given in \eqref{R_linear}, that, when inserted in  Eq.~\eqref{FP_Final} yields
%
%\begin{align}  
%\label{FP_Final_OU} 
%&\partial_t {P}(x;t)=
%\left\{\gamma \partial_x x\, +D_f\partial _x^2 \right\}P(x;t)\nonumber \\
%&+ \epsilon^2 \partial_x x
%\left\{
%\partial_x\,x
%+D
% \frac{2 \tau}{2 \gamma  \tau +1}
%\partial_x^2
%\right\} P(x;t).
%\end{align}  
%%
 %
\section{Conclusions\label{sec:conclusion}}
The Fokker-Planck Equation (FPE) holds a central position in statistical mechanics. Initially derived as the Kramers-Moyal expansion of the 
Master Equation (ME), limited to Markovian systems, it's recognized as applicable to non-Markovian processes. 
Indeed, the FPE emerges by eliminating irrelevant or fast variables, weakly interacting with the part of interest, through perturbation 
techniques like Zwanzig and Mori's projective methods, or considerations on the order of magnitude of the generalized cumulants.  
Thus, it stands as the most important equation to derive the PDF 
time evolution in these approximations.

Moreover, the FPE has the advantage of being a second order classical parabolic PDE with well-studied properties on solution existence and positivity. 
Its importance and widespread use are undeniable.

The connection between FPE and Fick's law is not coincidental. The FPE, when expressed as a continuity equation, reveals that the current linked to the stochastic process involves a diffusion term that is proportional to the 
gradient of the PDF, constituting Fick's law. Conversely, assuming Fick's law holds, the continuity equation emerges as a second-order PDE, exhibiting the structure of an FPE. 
Hence, whether Fick's law holds or not, and the ME with FPE structure structure are intricately connected.

The extensive use of the FPE has led to the development of numerous methods for extracting 
crucial statistical information. Standard spectral analysis procedures, similar to those 
applied to the Schr\"odinger equation in quantum mechanics, can be employed. Additionally, the 
diffusion and drift coefficients of the FPE allow for the derivation of an analytical 
expression for the mean first-passage time. This important quantity represents the 
average time for a 
trajectory, starting from an initial position $x_0$, to reach a specified target point $x_T$ 
for the first time.

In our study, we demonstrated that when the system of interest is inherently noisy-featuring 
sources like Nyquist noise in electric circuits, 
various thermal fluctuations, rapid internal dynamics, intrinsic measurement errors, etc., the
mentioned standard procedures for eliminating fast or 
weakly interacting variables (often, but not necessarily,  modeled as stochastic processes) 
lead to a third-order PDE, instead of an FPE. 
In fact, an additive third-order partial differential operator emerges from the interplay 
between the standard diffusion process due to internal noise and the  
diffusion process due to the external colored stochastic process (or to the irrelevant degrees 
of freedom we project out).

Given the inevitability of such internal noise (of varying intensity), we conclude that 
the third-order PDE 
should be considered more fundamental than the FPE in statistical physics. This fact also implies 
the breakdown of the Fick's law.

While this approach can be extended to accommodate more general drift fields, our current focus in this work is on the simpler linear drift case, which is widely employed across various disciplines.  
The analytical expressions of the moments of the PDF reveals 
that the unexpected third derivative term
 significantly tightens the 
 equilibrium PDF, in comparison to what we would obtain if we dropped this term, maintaining just the standard structure of the FPE.
 Figures \ref{fig:PDF_OU1} and \ref{fig:PDF_OU2} support this observation, showing perfect 
 agreement between numerical simulations of 
 the SDE and the third-order PDE.
In particular, the figures clearly illustrate that the actual PDFs, effectively 
captured by the third-order PDE, exhibit a significantly more narrow equilibrium
PDFs compared to those derived from the FPE. 
This observation agrees with the general finding, emphasized in section~\ref{sec:third_order}, 
that the third-order differential contribution to ME leads to a reduction in the moments of $x$.
Consequently, the tails of the actual equilibrium PDF (and those of the equilibrium PDF
of the third-order PDE) decay more rapidly than those of the FPE, indicating that crucial 
statistical quantities, such as the average first-passage time, computed using standard FPE 
techniques, would yield inaccurate results.

Thus, the fact that for the statistical behavior of a specific part of a complex system, 
the third-order PDE should be considered 
more fundamental than the FPE raises the question of how to extend to this PDE
the general methods and 
results, such like those which allow the derivation 
of relevant statistical information from the FPE. 
For example, in the 1-D case, it would be interesting to obtain an analytical expression for 
the equilibrium PDF or a closed expression for the mean first-passage time. 
All these are matter of future works.
\section{Aknowledgement} 
We thank the Green Data Center of University of Pisa for providing the computational power needed for the present paper. This research work was supported in part by ISMAR-CNR and UniPi
institutional funds. M.B. acknowledges financial support from UTA Mayor
Project No. 8738-23.
\appendix
\section{The multidimensional case\label{app:multiD}}
In this Appendix we generalize the result \eqref{FP_Final} to the multi-dimensional case. 
For the reader convenience, we reapt here the general  $N$-D extension of the SDE \eqref{SDE2_1}, already introduced in \eqref{SDE2_MD_}:
\begin{align}
\label{SDE2_MD}
\dot {\vect x}= -\mathbb{E}\cdot \vect x+\vect f(t)
-    \bm{\Xi}(t)\cdot \vect x
\end{align}
where $\bm x:=(x_1,\dots, x_N$), $\mathbb{E}$ and $\bm{\Xi}(t)$ are $N\times N$ matrices with constant  %$\mathbb{E}_{i,j}%=e_{i,j}$ 
and stochastic components% 
%$\bm{\Xi}_{i,j}(t)=\xi_{i,j}(t)$
, respectively%
%and $i,j\in\{1,2,...,N\}$
. Moreover, $\vect f(t)$ is a multidimensional white noise with correlation, or diffusion matrix given by $\mathbb{D}$.

As for the 1-D case, to any realization
$\bm \Xi(\cdot)$ of the matrix stochastic process $\bm \Xi(u)$, $0\le u \le t$, 
from  \eqref{SDE2_MD} we can write 
  the following  Liouville equation  for the PDF of $\vect x$, that we indicate with 
  ${P}_{\bm \Xi(\cdot)}(\vect x,t)$:
\begin{align}
\label{Liouville_MD}
&\partial_t {P}_{\bm \Xi(\cdot)}(\vect x,t) \nonumber \\
&= 
\{{\cal L}_a+  {\cal  L}_{\bm \Xi(t)}\}  
{P}_{\bm \Xi(\cdot)}(\vect x,t) ,
\end{align}
in which the unperturbed Liouvillian is
($\vect \partial$ is the $N$-D gradient operator and
%and repeated indices implies summation from $1$ to $N$, moreover, 
the superscript ``$T$'' means ``transpose''):
\begin{equation}
\label{La_MD}
{\cal L}_a:=%\partial_i e_{i,k} x_k+\partial_i D_{i,k} \partial_k=
\vect \partial^T\cdot \mathbb{E} \cdot \vect x+
+\vect \partial^T\cdot \mathbb{D} \cdot \vect \partial
\end{equation}
and the Liouville perturbation operator is
\begin{equation}
\label{LI_MD}
{\cal L}_{\bm \Xi(t)}:=%\partial_i \xi_{i,k}(t)  x_k=
\vect \partial^T\cdot \bm{\Xi}(t) \cdot \vect x.
\end{equation}
We rewrite the Liouville equation \eqref{Liouville_MD} 
 in interaction representation:
\begin{equation}  
        \partial_t      \tilde {P}_{\bm \Xi(\cdot)}(\vect x,t) 
        = \tilde{\cal L}_{\bm \Xi(t)}(t) 
        \tilde {P}_{\bm \Xi(\cdot)}(\vect x,t)  ,  
        \label{app:liouvilleequation1_MD}  
\end{equation}  
where 
\begin{equation} 
        \tilde {P}_{\bm \Xi(\cdot)}(\vect x,t)  :=  
        e^{-{\cal L}_a t} {P}_{\bm \Xi(\cdot)}(\vect x,t) 
\end{equation}
and
\begin{align}
\label{LIt_MD_0}
\tilde{\cal L}_{\bm \Xi(t)}(t) :=& 
e^{-{\cal L}_a t}{\cal L}_{\bm \Xi(t)} e^{{\cal L}_a t}
=e^{-{\cal L}_a^\times t} [{\cal L}_{\bm \Xi(t)}].
\end{align}
Integrating (\ref{app:liouvilleequation1_MD}) and averaging over the
realization of $\bm \Xi (t)$, we get
\begin{equation}
        \label{app:texp}
        \tilde {P}(\vect x;t) = 
        \langle \overleftarrow{\exp} \left[ \int_0^t du \;  \tilde{\cal L}_{\bm \Xi(t)} (u) \right]
        \rangle_{{\bm \Xi}}\,
        {P}(\vect x;0)
\end{equation}
in which $\overleftarrow{\exp} [...]$ is the standard chronological ordered 
exponential (from right to left) and
$       \tilde P(\vect x;t):=  
e^{-{\cal L}_a t} P(\vect x;t)$ with 
$P(\vect x;t):=
\langle{P}_{\bm \Xi(\cdot)}(\vect x,t)\rangle_{\bm \Xi}$.
By using the generalized cumulant approach and retaining only the second cumulant we get 
the following ME for the PDF of $\vect x$: 
\begin{equation}  
\label{FP_F2_MD} 
\partial_t {P}(\vect x;t)= {\cal L}_a P(x;t)+   
\int_0^{\infty}\text{d}u\,
\langle 
{\cal L}_{\bm \Xi(t)} \tilde{\cal L}_{\bm \Xi(-u)}(-u)
\rangle_{\bm \Xi}\;
{P} (\vect x;t),
\end{equation}  
corresponding to the $N$-D version of  Eq.~\eqref{FP_F2}.
To obtain the explicit expression, as PDE, of the ME \eqref{FP_F2_MD}, we must solve the Lie evolution of 
${\cal L}_{\bm \Xi(t)}$ along the Liouvillian ${\cal L}_a$, i.e. we have to explicitly evaluate 
$\tilde{\cal L}_{\bm \Xi(u)}(u)$ of \eqref{LIt_MD_0}, in which  ${\cal L}_a$ and 
${\cal L}_{\bm \Xi(t)}$ are given in \eqref{La_MD} and \eqref{LI_MD}, respectively. For that, let us 
start considering the operator identity 
$e^{({\cal L}_A+{\cal L}_B) t}=e^{{\cal L}_A t} \cdot \overleftarrow{\exp}\left(
\int_0^t \text{d} u\,\bar{\cal L}_B(u) \right)$ in which  
${\cal L}_A$ and ${\cal L}_B$ are operators that in general do not commute with each other, and where
$\bar{\cal L}_B(u):=e^{-{\cal L}_A u} {\cal L}_B e^{{\cal L}_Au}$.
From this identity, by making the
associations ${\cal L}_A=\vect \partial^T\cdot \mathbb{E} \cdot \vect x$
and ${\cal L}_B=\vect \partial^T\cdot \mathbb{D} \cdot \vect \partial$ (thus, 
${\cal L}_A+{\cal L}_B={\cal L}_a$),  with a few algebra we easily obtain:
%
%\begin{equation}
%\label{expLa_MD}
%e^{{\cal L}_a t}=
%e^{\vect \partial^T\cdot \mathbb{E} \cdot \vect x\,t}
%\cdot 
%\overleftarrow{\exp}\left(
%\int_0^t \text{d} u\,
%\vect \partial^T\cdot \mathbb{O}(u)\cdot \mathbb{D}\cdot \mathbb{O}^T(u) \cdot \vect \partial\, 
%\right)
%\end{equation}
% 
%
\begin{equation}
\label{expLa_MD}
e^{{\cal L}_a t}=
e^{\vect \partial^T\cdot \mathbb{E} \cdot \vect x\,t}
\cdot 
\overleftarrow{\exp}\left(
\int_0^t \text{d} u\,
\vect \partial^T\cdot e^{\mathbb{E}u}\cdot \mathbb{D}\cdot e^{\mathbb{E}^Tu}\cdot \vect \partial\, 
\right).
\end{equation}
% 
%where
%\begin{equation}
%\label{O}
%\mathbb{O}(u):=e^{\mathbb{E}u}.
%\end{equation}
%
By using \eqref{expLa_MD} and  \eqref{La_MD} in \eqref{LIt_MD_0} we get
\begin{align}
\label{LIt_MD}
\tilde{\cal L}_{\bm \Xi(t)}(t) :=& 
e^{-{\cal L}_a t}{\cal L}_{\bm \Xi(t)} e^{{\cal L}_a t}
\nonumber \\
=&e^{-\vect \partial^T\cdot \mathbb{E} \cdot \vect x\,t}
\cdot 
\overleftarrow{\exp}\left(-
\int_{-t}^0 \text{d} u\,
\vect \partial^T\cdot e^{\mathbb{E}u}\cdot \mathbb{D}\cdot e^{\mathbb{E}^Tu} \cdot \vect \partial\, 
\right)\nonumber\\
&\cdot
{\cal L}_{\bm \Xi(t)}\cdot
e^{\vect \partial^T\cdot \mathbb{E} \cdot \vect x\,t}
\cdot 
\overleftarrow{\exp}\left(
\int_0^t \text{d} u\,
\vect \partial^T\cdot e^{\mathbb{E}u}\cdot \mathbb{D}\cdot e^{\mathbb{E}^Tu} \cdot \vect \partial\, 
\right) \nonumber \\
=& 
\overleftarrow{\exp}\left(-
\int_{0}^t \text{d} u\,
\vect \partial^T\cdot e^{\mathbb{E}u}\cdot \mathbb{D}\cdot e^{\mathbb{E}^Tu} \cdot \vect \partial\, 
\right)^\times
\left[
e^{-\vect \partial^T\cdot \mathbb{E} \cdot \vect x\,t^\times}
\left[{\cal L}_{\bm \Xi(t)}\right]\right].
\end{align}
In the last side of the following equation we have exploited the following identity, easily demonstrated:
\begin{align}
&e^{-{\cal L}_A \theta} \cdot \overleftarrow{\exp}\left(-
\int_0^t \text{d} u\,\bar{\cal L}_B(u) \right)
=
e^{-{\cal L}_A \theta^\times} \left[ 
\overleftarrow{\exp}\left(-
\int_0^t \text{d} u\,\bar{\cal L}_B(u) \right)
\right]
e^{-{\cal L}_A \theta} 
\nonumber \\
&=
\overleftarrow{\exp}\left(-
\int_\theta^{t+\theta} \text{d} u\,\bar{\cal L}_B(u) \right)
e^{-{\cal L}_A \theta}.
\end{align}
By using \eqref{LI_MD} and the results of~\cite{bJMP59}, in particular those in Section~VIA, we have 
\begin{align}
e^{-\vect \partial^T\cdot \mathbb{E} \cdot \vect x\,t^\times}
\left[{\cal L}_{\bm \Xi(t)}\right]&=
e^{-\vect \partial^T\cdot \mathbb{E} \cdot \vect x\,t^\times}
\left[
\vect \partial^T\cdot \bm{\Xi}(t) \cdot \vect x
\right]%\nonumber \\
%&
=
%\vect \partial^T\cdot
%e^{\mathbb{E}t}\cdot
%\bm{\Xi}(t)\cdot
%e^{-\mathbb{E}t}\cdot \vect x= 
\vect \partial^T\cdot
e^{\mathbb{E}t^\times}
\left[\bm{\Xi}(t)\right]
\cdot \vect x.
\end{align}
Inserting this result in \eqref{LIt_MD} we obtain
\begin{align}
\label{LIt_MD_2}
&\tilde{\cal L}_{\bm \Xi(t)}(t) 
=\nonumber \\
&\overleftarrow{\exp}\left(-
\int_{0}^t \text{d} u\,
\vect \partial^T\cdot e^{\mathbb{E}u}\cdot \mathbb{D}\cdot e^{\mathbb{E}^Tu} \cdot \vect \partial\, 
\right)^\times
\left[
\vect \partial^T\cdot
e^{\mathbb{E}t^\times}
\left[\bm{\Xi}(t)\right]
\cdot \vect x\right].
\end{align}
By expanding the above series of nexted commutators we see that all the terms are zero, apart the zeroth and the first ones. Therefore,  we get
\begin{align}
\label{LIt_MD_3}
&\tilde{\cal L}_{\bm \Xi(t)}(t) 
=
\vect \partial^T\cdot
e^{\mathbb{E}t^\times}
\left[\bm{\Xi}(t)\right]
\cdot \vect x
\nonumber \\
&-\vect \partial^T\cdot
e^{\mathbb{E}t^\times}
\left[\bm{\Xi}(t)\right]
\cdot
\int_{0}^t \text{d} u\,\left\{
\left(
e^{\mathbb{E}u}\cdot \mathbb{D}\cdot e^{\mathbb{E}^Tu}
\right)^T
+e^{\mathbb{E}u}\cdot \mathbb{D}\cdot e^{\mathbb{E}^Tu}
\right\} \cdot 
\vect \partial,
\end{align}
that, given the symmetry property of the diffusion coefficient matrix, yields the final
explicit differential form for the interaction representation of the Liouvillian 
${\cal L}_{\bm \Xi(t)}$:
\begin{align}
\label{LIt_MD_4}
\tilde{\cal L}_{\bm \Xi(t)}(t) =
\vect \partial^T\cdot
e^{\mathbb{E}t^\times}
\left[\bm{\Xi}(t)\right]
\cdot
\left\{\vect x-2
\int_{0}^t \text{d} u\,
e^{\mathbb{E}u}\cdot \mathbb{D}\cdot e^{\mathbb{E}^Tu}
\cdot \vect \partial 
\right\}.
\end{align}
Thus, by using this expression in the ME \eqref{FP_F2_MD}, together with Eqs.~\eqref{La_MD} and \eqref{LI_MD},
we arrive to the final general PDE of third order for the PDF of $\vect x$ for the multi-dimensional case:
\begin{align}  
\label{FP_F2_MD_2} 
&\partial_t {P}(\vect x;t)= \left\{\vect \partial^T\cdot \mathbb{E} \cdot \vect x+
+\vect \partial^T\cdot \mathbb{D} \cdot \vect \partial\right\} P(\vect x;t)
+\int_0^{\infty}\text{d}u  \nonumber \\
&\times
\langle
\vect \partial^T\cdot \bm{\Xi}(t) \cdot \vect x 
\left(
\vect \partial^T\cdot
e^{-\mathbb{E}u^\times}
\left[\bm{\Xi}(-u)\right]
\cdot
\left\{\vect x+2
\int_{0}^u \text{d} u\,
e^{\mathbb{E}u}\cdot \mathbb{D}\cdot e^{\mathbb{E}^Tu}
\cdot \vect \partial 
\right\}
\right)\rangle_{\bm \Xi}
\;{P} (\vect x;t),
\end{align}  
In the simplified case in which $\bm \Xi(t)=\epsilon \mathbb{G}\xi(t)$, with 
$\langle \xi(t)\xi \rangle_{\bm \Xi}=\varphi(t)$ then we have
\begin{align}  
\label{FP_F2_MD_3} 
&\partial_t {P}(\vect x;t)= \left\{\vect \partial^T\cdot \mathbb{E} \cdot \vect x+
+\vect \partial^T\cdot \mathbb{D} \cdot \vect \partial\right\} P(x;t)
+\frac{\delta^2}{\tau^2} \int_0^{\infty}\text{d}u  \,\varphi(u)\nonumber \\
&\times
\vect \partial^T\cdot \mathbb{G} \cdot \vect x 
\left(
\vect \partial^T\cdot
e^{-\mathbb{E}u^\times}
\left[\mathbb{G}\right]
\cdot
\left\{\vect x+2
\int_{0}^u \text{d} u\,
e^{\mathbb{E}u}\cdot \mathbb{D}\cdot e^{\mathbb{E}^Tu}
\cdot \vect \partial 
\right\}
\right)
\;{P} (\vect x;t),
\end{align}  
where we have also used the definition of the adimensional  parameter $\delta:=\epsilon \tau$, that is the relevant small quantity in the cumulant expansion.
\section{The cumulant approach as a systematic way to obtain a ME for the reduced PDF of $x$\label{app:A}}

In this Appendix we  outline a few minima key steps to obtain the
FPE~\eqref{FPE_Linear} and the ME~\eqref{FP_F2}, starting from the generalized cumulant 
(or $M$-cumulant) approach formally presented in \citep{bbJSTAT4}. 
We begin with the generic Liouville equation \eqref{Liouville}
(the stochastic process is one-dimensional, but the extension to multi-dimensional cases is straightforward), expressed in interaction representation:
\begin{equation}
\partial_t \tilde {P}_{\xi(\cdot)}(x,t)
= \epsilon \xi (t),\tilde{\cal L}_I(t) \tilde {P}_{\xi(\cdot)}(x,t).
\label{app:liouvilleequation1}
\end{equation}
Here,
\begin{equation}
\tilde {P}_{\xi(\cdot)}(x,t) :=
e^{-{\cal L}_a t} {P}_{\xi(\cdot)}(x,t)
\end{equation}
and
\begin{equation}
\label{app:LIt}
\tilde{\cal L}_I(t) := e^{-{\cal L}_a t}\tilde{\cal L}_I e^{{\cal L}_a t}
=e^{-{\cal L}_a^\times t} [{\cal L}_I].
\end{equation}
In~\cite{bJMP59}, $\tilde{\cal L}_I(t)$ of (\ref{app:LIt}) is also referred to as the Lie evolution of the operator ${\cal L}_I$ along the Liouvillian 
${\cal L}_a$, for a time $-t$. 

Integrating (\ref{app:liouvilleequation1}) and averaging over the
realization of $\xi (t)$, we get
\begin{equation}
        \label{app:texp}
        \tilde P(x;t) = 
        \langle \overleftarrow{\exp} \left[\epsilon \int_0^t du \; \xi (u) \tilde{\cal L}_I (u) \right]
        \rangle_{\xi} P(x;0)
\end{equation}
in which $\overleftarrow{\exp} [...]$ is the standard chronological ordered 
exponential (from right to left) and
$       \tilde P(x;t):=  
e^{-{\cal L}_a t} P(x;t)$ with $P(x;t):=\langle{P}_{\xi(\cdot)}(x,t)\rangle_\xi$. Moreover,
 we have exploited the assumption  that at the initial time $t=0$ the total PDF
 factorizes as ${P}_{\xi(\cdot)}(x,0)=P(x;0)p(\xi)$. 
 This is equivalent to stating that at the initial time the PDF of $x$ does not depend on the 
 possible values of the process $\xi$, 
 or alternatively,  we wait long enough so that the initial conditions became irrelevant. 
Apart that, Eq.~(\ref{app:texp}) is exact; no approximations have been introduced at this level.

We can look at the r.h.s. of \eqref{app:texp} as a 
sort of characteristic function, or moment generating function, with wave number $k:=\text{i}\epsilon$, for the stochastic operator
\begin{equation}
    \label{Omega}
\Omega(u):=\xi (u) \tilde{\cal L}_I(u).
\end{equation}
Formally, we can then introduce a generalized
cumulant generating  function~\citep{bbJSTAT4}:
\begin{equation}
        \label{app:texp2}
        \langle \overleftarrow{\exp} \left[\epsilon \int_0^t du \; \xi (u) \tilde{\cal L}_I (u) \right]
        \rangle_{\xi} := 
        \overleftarrow{\exp} \left[{\cal K}(\epsilon,t) \right]
\end{equation}
with 
\begin{equation}
        \label{K_expansion}
        {\cal K}(\epsilon,t)=\sum_{i=1}^\infty\epsilon^i{\cal K}_i(t).
\end{equation}
 As for standard stochastic processes, we define the
$n$-times joint $M$-cumulant of $\Omega(u)$, that we indicate as 
$\langle\langle \Omega(u_1) \Omega(u_2)...\Omega(u_n)  \rangle\rangle$, by setting
\begin{align}
        \label{k_int}
        {\cal K}_i (t):=& \int_0^t \text{d}u_1 \int_{0}^{u_1} \text{d}u_2
        ... \int_{0}^{u_{n-1}} \text{d}u_n
        \langle\langle \Omega(u_1) \Omega(u_2)...\Omega(u_n)  \rangle\rangle.
\end{align}
Using  \eqref{k_int} in the r.h.s. of \eqref{app:texp2} and expanding both exponential functions, 
we get the standard relationship among cumulants and moments.
For example, the joint two and four times   $M$-cumulants are given in terms of moments as
(to improve readability, until the end of this paragraph we will avoid putting the subscript ``$\xi$' to the angle brackets):
\begin{align}
        \label{two}
         \langle\langle\Omega(u_1)\,\Omega(u_2)\rangle\rangle=  \langle\Omega(u_1)\,\Omega(u_2)\rangle=
         \tilde{\cal L}_I(u_1)\tilde{\cal L}_I(u_2)
\langle   \xi (u_1)\,\xi (u_2)\rangle
\end{align}
and
\begin{align}
        \label{four}
        &\langle\langle\Omega(u_1)\,\Omega(u_2)\Omega(u_3)\,\Omega(u_4)\rangle\rangle=\nonumber \\
        &
        \langle\Omega(u_1)\,\Omega(u_2) \Omega(u_3)\,\Omega(u_4)\rangle
        -       \langle\Omega(u_1)\,\Omega(u_2)\rangle \langle\Omega(u_3)\,\Omega(u_4)\rangle\nonumber \\
        &-      \langle\Omega(u_1)\,\Omega(u_3)\rangle \langle\Omega(u_2)\,\Omega(u_4)\rangle
        -       \langle\Omega(u_1)\,\Omega(u_4)\rangle \langle\Omega(u_2)\,\Omega(u_3)\rangle=\nonumber \\
        &       
        \, \tilde{\cal L}_I(u_1)\tilde{\cal L}_I(u_2)\tilde{\cal L}_I(u_3)\tilde{\cal L}_I(u_4)
                \left[\langle   \xi (u_1)\,\xi (u_2)\,\xi (u_3)\,\xi (u_4)\rangle-\langle \xi (u_1)\xi (u_2)
        \langle\xi (u_3)\xi (u_4)\rangle \right]\nonumber \\
        &-\tilde{\cal L}_I(u_1)\tilde{\cal L}_I(u_3)\tilde{\cal L}_I(u_2)\tilde{\cal L}_I(u_4)
        \langle\xi (u_1)\xi (u_3)\rangle  \langle\xi (u_2)\xi (u_4) \rangle\nonumber\\ 
        &-\tilde{\cal L}_I(u_1)\tilde{\cal L}_I(u_4)\tilde{\cal L}_I(u_2)\tilde{\cal L}_I(u_3)
        \langle\xi (u_1)\xi (u_4)\rangle  \langle\xi (u_2)\xi (u_3)\rangle,
\end{align}
respectively.
From \eqref{four} it is clear that the Gaussian nature of $\xi(t)$ does not 
implies the same for $\Omega(t)$ of Eq.~\eqref{Omega}, as the time-dependent Liouvillian 
$\tilde{\cal L}_I(u)$ generally does not commute with itself evaluated at different times. 
However, when the unperturbed Liouvillian ${\cal L}_a$ and perturbation Liouvillian 
${\cal L}_I$ commute with each other, as in the case of Eq.~\eqref{La} with $D_f =0$ and ${\cal L}_I$ of 
Eq.~\eqref{LI}, we have $\tilde{\cal L}_I(u)={\cal L}_I$, that does not depend on time. 
Hence, in this case the Gaussian nature of $\xi(t)$ is transferred to the stochastic 
operator $\Omega(t)$. Therefore,
in this specific scenario, the $M$-cumulant series appearing in the exponential function 
of \eqref{app:texp2} reduces to only the second term containing the second $M$-cumulant, 
simplifying to (without loss of generality, we consider the average value of 
$\xi (t)$ to be zero):
 \begin{equation}
        \label{app:texp_G}
        \tilde P(x;t) = 
 {\exp} \left[\epsilon^2  {\cal L}_I {\cal L}_I
 \int_0^t \text{d}u_1 \int_{0}^{u_1} \text{d}u_2
\langle   \xi (u_1)\,\xi (u_2)\rangle
\right] P(x;0).
 \end{equation}
Time-deriving this result we obtain
\begin{align}  
        \partial_t \tilde P(x;t)
        =&
        \epsilon^2 {\cal L}_I {\cal L}_I
        \int_0^t  \text{d}u
                \langle   \xi (t)\,\xi (u)\rangle\;
        \tilde P(x;t)\nonumber \\
        =&
        \epsilon^2 {\cal L}_I {\cal L}_I
\tau    \tilde P(x;t).
        \label{app:liouvilleequation_G}  
\end{align}  
Getting rid of the interaction representation and by using \eqref{La} with $D_f =0$
and  \eqref{LI}, Eq.~\eqref{app:liouvilleequation_G} becomes exactly the 
FPE~\eqref{FPE_Linear}.

In the more general case, the Liouvillians ${\cal L}_a$ and ${\cal L}_I$ do not commute 
with each other, so $\tilde{\cal L}_I(u)$ of \eqref{app:LIt} depends on time. 
The advantage of utilizing the $M$-cumulants lies in the fact that, similar to 
standard cumulants, they are exactly zero when referring to independent random 
variables~\citep{bbJSTAT4}. 
Thus, the time lag between two events increases until they become independent of each 
other, any joint $M$-cumulant containing these two events must tend to zero. To model 
this situation more realistically, we assume that independence does not occur abruptly 
at a fixed time lag $\bar \tau$ but instead follows a smoother pattern, characterized by
an exponential trend. Formally, for a series of events 
$\xi(t_1), \xi(t_2),...,\xi(t_n)$
with $t_1\ge t_2\ge...\ge t_n$, we assume that the corresponding joint $n$-cumulant 
decays at least exponentially with the time lag $u_1-u_n$::
\begin{equation}
        \label{exp-correlated}
\left|\langle\langle \Omega(u_1) \Omega(u_2)...\Omega(u_n)  \rangle\rangle\right|
\lesssim \exp(-(u_1-u_n)/\bar \tau).
\end{equation}
In this scenario, along with the 
definitions \eqref{K_expansion} and \eqref{k_int}, it is evident that the argument of 
the exponential function in the right-hand side of \eqref{app:texp2} now yields a power 
series of  $\bar \delta:=\epsilon\bar \tau$.
For a sufficiently small $\bar \delta$, we can truncate this series to the first 
non-zero term, which is the second one. 
Thus, Eq.~\eqref{app:texp2}, combined with Eq.~\eqref{app:texp2} and \eqref{two}, gives
 \begin{equation}
        \label{app:texp_}
        \tilde P(x;t) = 
 \overleftarrow{\exp} \left[\epsilon^2
 \int_0^t \text{d}u_1 \int_{0}^{u_1} \text{d}u_2
 \tilde{\cal L}_I(u_1)\tilde{\cal L}_I(u_2)
\langle   \xi (u_1)\,\xi (u_2)\rangle+ O(\bar \delta^4)
\right] P(x;0).
 \end{equation}
Time-deriving this result we obtain
\begin{equation}  
        \partial_t \tilde P(x;t)
        =\epsilon^2
        \int_0^t  \text{d}u
        \tilde{\cal L}_I(t)\tilde{\cal L}_I(u)
        \langle   \xi (t)\,\xi (u)\rangle\;
        \tilde P(x;t)+ O(\bar \delta^4 t/\bar \tau)
        \label{app:liouvilleequation2}  
\end{equation}  
Getting rid of the interaction representation and by using again \eqref{La} (but now letting   $D_f\ne0$)
and  \eqref{LI}, Eq.~\eqref{app:liouvilleequation2} becomes  the approximate
ME~\eqref{FP_F2}.
%

%While general results for evaluating the Lie evolution of differential operators can be found in~\cite{bJMP59}, 
%for the present simple 1-D case, it is straightforward to obtain:\footnote{
%%
%%BEGIN FOOTNOTE
%From the definitions of \eqref{La} and \eqref{LI}, and %by taking into account that the Lie evolution in the %l.h.s. of 
%\eqref{app:LIt2} is a power series of commutators, we have the following trivial series of equations: 
%$e^{-{\cal L}_a^\times t} [{\cal L}_I] =e^{-{\cal L}_a^\times t} [\partial_x ]e^{-{\cal L}_a^\times t} [I(x)]
%       =e^{-{\cal L}_a^\times t} \left[\partial_x \frac{C(x)}{C(x)}\right]I(x_0(x,u))
%       = \partial_x C(x)\, e^{-{\cal L}_a^\times t}\left[\frac{1}{C(x)}\right]I(x_0(x,u))
%       =\partial_x  \,\frac{ C(x)}{C(x_0(x,u))}I(x_0(x,u))$
%}
%%
%%END FOOTNOTE}:
%%
%\begin{align}
%       \label{app:LIt2} 
%       e^{-{\cal L}_a^\times t} [{\cal L}_I] 
%       =\partial_x  \,\frac{ C(x)}{C(x_0(x,u))}I(x_0(x,u)),
%\end{align}
%
%where $x_0(x;-u)$ is the backwards unperturbed trajectory evolved for a time $u$, starting from
%the initial position $x$.

%
%
%
%
\newpage
%+Bibliography
%\begin{thebibliography}{99}
%\bibitem{Label1} ...
%\bibitem{Label2} ...
%\end{thebibliography}
%-Bibliography
%
%\bibliography{BiblioCentrale_}% Produces the bibliography via BibTeX.

\begin{thebibliography}{32}%
\makeatletter
\providecommand \@ifxundefined [1]{%
 \@ifx{#1\undefined}
}%
\providecommand \@ifnum [1]{%
 \ifnum #1\expandafter \@firstoftwo
 \else \expandafter \@secondoftwo
 \fi
}%
\providecommand \@ifx [1]{%
 \ifx #1\expandafter \@firstoftwo
 \else \expandafter \@secondoftwo
 \fi
}%
\providecommand \natexlab [1]{#1}%
\providecommand \enquote  [1]{``#1''}%
\providecommand \bibnamefont  [1]{#1}%
\providecommand \bibfnamefont [1]{#1}%
\providecommand \citenamefont [1]{#1}%
\providecommand \href@noop [0]{\@secondoftwo}%
\providecommand \href [0]{\begingroup \@sanitize@url \@href}%
\providecommand \@href[1]{\@@startlink{#1}\@@href}%
\providecommand \@@href[1]{\endgroup#1\@@endlink}%
\providecommand \@sanitize@url [0]{\catcode `\\12\catcode `\$12\catcode
  `\&12\catcode `\#12\catcode `\^12\catcode `\_12\catcode `\%12\relax}%
\providecommand \@@startlink[1]{}%
\providecommand \@@endlink[0]{}%
\providecommand \url  [0]{\begingroup\@sanitize@url \@url }%
\providecommand \@url [1]{\endgroup\@href {#1}{\urlprefix }}%
\providecommand \urlprefix  [0]{URL }%
\providecommand \Eprint [0]{\href }%
\providecommand \doibase [0]{http://dx.doi.org/}%
\providecommand \selectlanguage [0]{\@gobble}%
\providecommand \bibinfo  [0]{\@secondoftwo}%
\providecommand \bibfield  [0]{\@secondoftwo}%
\providecommand \translation [1]{[#1]}%
\providecommand \BibitemOpen [0]{}%
\providecommand \bibitemStop [0]{}%
\providecommand \bibitemNoStop [0]{.\EOS\space}%
\providecommand \EOS [0]{\spacefactor3000\relax}%
\providecommand \BibitemShut  [1]{\csname bibitem#1\endcsname}%
\let\auto@bib@innerbib\@empty
%</preamble>
\bibitem [{\citenamefont {Deutsch}(1994{\natexlab{a}})}]{dPA208b}%
  \BibitemOpen
  \bibfield  {author} {\bibinfo {author} {\bibfnamefont {J.}~\bibnamefont
  {Deutsch}},\ }\bibfield  {title} {\enquote {\bibinfo {title} {Probability
  distributions for multicomponent systems with multiplicative noise},}\ }\href
  {\doibase https://doi.org/10.1016/0378-4371(94)00054-9} {\bibfield  {journal}
  {\bibinfo  {journal} {Physica A: Statistical Mechanics and its Applications}\
  }\textbf {\bibinfo {volume} {208}},\ \bibinfo {pages} {445--461} (\bibinfo
  {year} {1994}{\natexlab{a}})}\BibitemShut {NoStop}%
\bibitem [{\citenamefont {Nakao}(1998)}]{nPRE58}%
  \BibitemOpen
  \bibfield  {author} {\bibinfo {author} {\bibfnamefont {H.}~\bibnamefont
  {Nakao}},\ }\bibfield  {title} {\enquote {\bibinfo {title} {{Asymptotic power
  law of moments in a random multiplicative process with weak additive
  noise}},}\ }\href {\doibase 10.1103/PhysRevE.58.1591} {\bibfield  {journal}
  {\bibinfo  {journal} {Physical Review E}\ }\textbf {\bibinfo {volume} {58}}
  (\bibinfo {year} {1998}),\ 10.1103/PhysRevE.58.1591}\BibitemShut {NoStop}%
\bibitem [{\citenamefont {Fujisaka}\ and\ \citenamefont
  {Yamada}(1985)}]{fyPTP74}%
  \BibitemOpen
  \bibfield  {author} {\bibinfo {author} {\bibfnamefont {H.}~\bibnamefont
  {Fujisaka}}\ and\ \bibinfo {author} {\bibfnamefont {T.}~\bibnamefont
  {Yamada}},\ }\bibfield  {title} {\enquote {\bibinfo {title} {{A New
  Intermittency in Coupled Dynamical Systems}},}\ }\href {\doibase
  10.1143/PTP.74.918} {\bibfield  {journal} {\bibinfo  {journal} {Progress of
  Theoretical Physics}\ }\textbf {\bibinfo {volume} {74}},\ \bibinfo {pages}
  {918--921} (\bibinfo {year} {1985})},\ \Eprint
  {http://arxiv.org/abs/https://academic.oup.com/ptp/article-pdf/74/4/918/5362%
597/74-4-918.pdf}
  {https://academic.oup.com/ptp/article-pdf/74/4/918/5362597/74-4-918.pdf}
  \BibitemShut {NoStop}%
\bibitem [{\citenamefont {Yamada}\ and\ \citenamefont
  {Fujisaka}(1987)}]{yfPLA8}%
  \BibitemOpen
  \bibfield  {author} {\bibinfo {author} {\bibfnamefont {T.}~\bibnamefont
  {Yamada}}\ and\ \bibinfo {author} {\bibfnamefont {H.}~\bibnamefont
  {Fujisaka}},\ }\bibfield  {title} {\enquote {\bibinfo {title} {Effect of
  inhomogeneity on intermittent chaos in a coupled system},}\ }\href {\doibase
  https://doi.org/10.1016/0375-9601(87)90545-7} {\bibfield  {journal} {\bibinfo
   {journal} {Physics Letters A}\ }\textbf {\bibinfo {volume} {124}},\ \bibinfo
  {pages} {421--425} (\bibinfo {year} {1987})}\BibitemShut {NoStop}%
\bibitem [{\citenamefont {Pikovsky}(1992)}]{Pikovsky1992}%
  \BibitemOpen
  \bibfield  {author} {\bibinfo {author} {\bibfnamefont {A.~S.}\ \bibnamefont
  {Pikovsky}},\ }\bibfield  {title} {\enquote {\bibinfo {title} {Statistics of
  trajectory separation in noisy dynamical systems},}\ }\href {\doibase
  10.1016/0375-9601(92)91049-W} {\bibfield  {journal} {\bibinfo  {journal}
  {Physics Letters A}\ }\textbf {\bibinfo {volume} {165}} (\bibinfo {year}
  {1992}),\ 10.1016/0375-9601(92)91049-W}\BibitemShut {NoStop}%
\bibitem [{\citenamefont {Platt}, \citenamefont {Hammel},\ and\ \citenamefont
  {Heagy}(1994)}]{phhPRL72}%
  \BibitemOpen
  \bibfield  {author} {\bibinfo {author} {\bibfnamefont {N.}~\bibnamefont
  {Platt}}, \bibinfo {author} {\bibfnamefont {S.~M.}\ \bibnamefont {Hammel}}, \
  and\ \bibinfo {author} {\bibfnamefont {J.~F.}\ \bibnamefont {Heagy}},\
  }\bibfield  {title} {\enquote {\bibinfo {title} {Effects of additive noise on
  on-off intermittency},}\ }\href {\doibase 10.1103/PhysRevLett.72.3498}
  {\bibfield  {journal} {\bibinfo  {journal} {Phys. Rev. Lett.}\ }\textbf
  {\bibinfo {volume} {72}},\ \bibinfo {pages} {3498--3501} (\bibinfo {year}
  {1994})}\BibitemShut {NoStop}%
\bibitem [{\citenamefont {Sornette}(1998)}]{sPRE57}%
  \BibitemOpen
  \bibfield  {author} {\bibinfo {author} {\bibfnamefont {D.}~\bibnamefont
  {Sornette}},\ }\bibfield  {title} {\enquote {\bibinfo {title}
  {{Multiplicative processes and power laws}},}\ }\href {\doibase
  10.1103/PhysRevE.57.4811} {\bibfield  {journal} {\bibinfo  {journal}
  {Physical Review E}\ }\textbf {\bibinfo {volume} {57}},\ \bibinfo {pages}
  {4811--4813} (\bibinfo {year} {1998})},\ \Eprint
  {http://arxiv.org/abs/9708231} {arXiv:9708231} \BibitemShut {NoStop}%
\bibitem [{\citenamefont {Sornette}\ and\ \citenamefont {Cont}(1997)}]{scJPII}%
  \BibitemOpen
  \bibfield  {author} {\bibinfo {author} {\bibfnamefont {D.}~\bibnamefont
  {Sornette}}\ and\ \bibinfo {author} {\bibfnamefont {R.}~\bibnamefont
  {Cont}},\ }\bibfield  {title} {\enquote {\bibinfo {title} {{Convergent
  multiplicative processes repelled from zero: Power laws and truncated power
  laws}},}\ }\href {\doibase 10.1051/jp1:1997169} {\bibfield  {journal}
  {\bibinfo  {journal} {Journal de Physique II}\ }\textbf {\bibinfo {volume}
  {7}} (\bibinfo {year} {1997}),\ 10.1051/jp1:1997169}\BibitemShut {NoStop}%
\bibitem [{\citenamefont {Schenzle}\ and\ \citenamefont
  {Brand}(1979)}]{sbPRA20}%
  \BibitemOpen
  \bibfield  {author} {\bibinfo {author} {\bibfnamefont {A.}~\bibnamefont
  {Schenzle}}\ and\ \bibinfo {author} {\bibfnamefont {H.}~\bibnamefont
  {Brand}},\ }\bibfield  {title} {\enquote {\bibinfo {title} {Multiplicative
  stochastic processes in statistical physics},}\ }\href {\doibase
  10.1103/PhysRevA.20.1628} {\bibfield  {journal} {\bibinfo  {journal} {Phys.
  Rev. A}\ }\textbf {\bibinfo {volume} {20}},\ \bibinfo {pages} {1628--1647}
  (\bibinfo {year} {1979})}\BibitemShut {NoStop}%
\bibitem [{\citenamefont {Graham}, \citenamefont {H\"ohnerbach},\ and\
  \citenamefont {Schenzle}(1982)}]{ghsPRL48}%
  \BibitemOpen
  \bibfield  {author} {\bibinfo {author} {\bibfnamefont {R.}~\bibnamefont
  {Graham}}, \bibinfo {author} {\bibfnamefont {M.}~\bibnamefont
  {H\"ohnerbach}}, \ and\ \bibinfo {author} {\bibfnamefont {A.}~\bibnamefont
  {Schenzle}},\ }\bibfield  {title} {\enquote {\bibinfo {title} {Statistical
  properties of light from a dye laser},}\ }\href {\doibase
  10.1103/PhysRevLett.48.1396} {\bibfield  {journal} {\bibinfo  {journal}
  {Phys. Rev. Lett.}\ }\textbf {\bibinfo {volume} {48}},\ \bibinfo {pages}
  {1396--1399} (\bibinfo {year} {1982})}\BibitemShut {NoStop}%
\bibitem [{\citenamefont {Levy}\ and\ \citenamefont
  {Solomon}(1996)}]{lsIJMPC07}%
  \BibitemOpen
  \bibfield  {author} {\bibinfo {author} {\bibfnamefont {M.}~\bibnamefont
  {Levy}}\ and\ \bibinfo {author} {\bibfnamefont {S.}~\bibnamefont {Solomon}},\
  }\bibfield  {title} {\enquote {\bibinfo {title} {Power laws are logarithmic
  "boltzmann" laws},}\ }\href {\doibase 10.1142/S0129183196000491} {\bibfield
  {journal} {\bibinfo  {journal} {International Journal of Modern Physics C}\
  }\textbf {\bibinfo {volume} {07}},\ \bibinfo {pages} {595--601} (\bibinfo
  {year} {1996})},\ \Eprint
  {http://arxiv.org/abs/https://doi.org/10.1142/S0129183196000491}
  {https://doi.org/10.1142/S0129183196000491} \BibitemShut {NoStop}%
\bibitem [{\citenamefont {Takayasu}, \citenamefont {Sato},\ and\ \citenamefont
  {Takayasu}(1997)}]{tstPRL79}%
  \BibitemOpen
  \bibfield  {author} {\bibinfo {author} {\bibfnamefont {H.}~\bibnamefont
  {Takayasu}}, \bibinfo {author} {\bibfnamefont {A.-H.}\ \bibnamefont {Sato}},
  \ and\ \bibinfo {author} {\bibfnamefont {M.}~\bibnamefont {Takayasu}},\
  }\bibfield  {title} {\enquote {\bibinfo {title} {Stable infinite variance
  fluctuations in randomly amplified langevin systems},}\ }\href {\doibase
  10.1103/PhysRevLett.79.966} {\bibfield  {journal} {\bibinfo  {journal} {Phys.
  Rev. Lett.}\ }\textbf {\bibinfo {volume} {79}},\ \bibinfo {pages} {966--969}
  (\bibinfo {year} {1997})}\BibitemShut {NoStop}%
\bibitem [{\citenamefont {Turelli}(1977)}]{tTPB12}%
  \BibitemOpen
  \bibfield  {author} {\bibinfo {author} {\bibfnamefont {M.}~\bibnamefont
  {Turelli}},\ }\bibfield  {title} {\enquote {\bibinfo {title} {Random
  environments and stochastic calculus},}\ }\href {\doibase
  https://doi.org/10.1016/0040-5809(77)90040-5} {\bibfield  {journal} {\bibinfo
   {journal} {Theoretical Population Biology}\ }\textbf {\bibinfo {volume}
  {12}},\ \bibinfo {pages} {140--178} (\bibinfo {year} {1977})}\BibitemShut
  {NoStop}%
\bibitem [{\citenamefont {Deutsch}(1994{\natexlab{b}})}]{dPA208}%
  \BibitemOpen
  \bibfield  {author} {\bibinfo {author} {\bibfnamefont {J.}~\bibnamefont
  {Deutsch}},\ }\bibfield  {title} {\enquote {\bibinfo {title} {Probability
  distributions for one component equations with multiplicative noise},}\
  }\href {\doibase https://doi.org/10.1016/0378-4371(94)00055-7} {\bibfield
  {journal} {\bibinfo  {journal} {Physica A: Statistical Mechanics and its
  Applications}\ }\textbf {\bibinfo {volume} {208}},\ \bibinfo {pages}
  {433--444} (\bibinfo {year} {1994}{\natexlab{b}})}\BibitemShut {NoStop}%
\bibitem [{\citenamefont {Lepri}(2020)}]{lCSF139}%
  \BibitemOpen
  \bibfield  {author} {\bibinfo {author} {\bibfnamefont {S.}~\bibnamefont
  {Lepri}},\ }\bibfield  {title} {\enquote {\bibinfo {title} {Chaotic
  fluctuations in graphs with amplification},}\ }\href {\doibase
  https://doi.org/10.1016/j.chaos.2020.110003} {\bibfield  {journal} {\bibinfo
  {journal} {Chaos, Solitons \& Fractals}\ }\textbf {\bibinfo {volume} {139}},\
  \bibinfo {pages} {110003} (\bibinfo {year} {2020})}\BibitemShut {NoStop}%
%\bibitem [{Note1()}]{Note1}%
\bibitem [{Note1()}]{Note1}%
  \BibitemOpen
  \bibinfo {note} {The general prescription is that there is a time $\protect
  \mathaccentV {bar}016\tau $ such that, for any time $t$, the instances of
  $\xi $ at times $t'>t+\protect \mathaccentV {bar}016\tau $ are ``almost
  statistically uncorrelated'' with the instances of $\xi $ at times $t'<t$.
  For ``almost statistically uncorrelated'' we mean that the joint probability
  density functions factorize up to terms $O(\protect \mathaccentV {bar}016\tau
  )$: $p_n(\xi ,t_1';\xi _2,t_2';...;\xi _k,t_k';\xi _{k+1},t_{1};...;\xi
  _n,t_h)= p_k(\xi ,t_1';\xi _2,t_2';...;\xi _k,t_k')\protect \tmspace
  +\thinmuskip {.1667em}p_h(\xi _{k+1},t_{1};...;\xi _n,t_h)+O(\protect
  \mathaccentV {bar}016{\tau })$ with $k,h,n\in \protect \mathbb {N}$, $k+h=n$
  and $t_i'> t_j+\protect \mathaccentV {bar}016{\tau }$. For example, $p_2(\xi
  ,t';\xi _2,t)=p_1(\xi ,t')\protect \tmspace +\thinmuskip {.1667em}p_1(\xi
  _2,t)+O(\protect \mathaccentV {bar}016{\tau })$.\label {note1}}\BibitemShut
  {NoStop}%
\bibitem [{\citenamefont {Lubashevsky}\ \emph {et~al.}(2010)\citenamefont
  {Lubashevsky}, \citenamefont {Heuer}, \citenamefont {Friedrich},\ and\
  \citenamefont {Usmanov}}]{lhuEPJB78}%
  \BibitemOpen
  \bibfield  {author} {\bibinfo {author} {\bibfnamefont {I.~A.}\ \bibnamefont
  {Lubashevsky}}, \bibinfo {author} {\bibfnamefont {A.}~\bibnamefont {Heuer}},
  \bibinfo {author} {\bibfnamefont {R.}~\bibnamefont {Friedrich}}, \ and\
  \bibinfo {author} {\bibfnamefont {R.}~\bibnamefont {Usmanov}},\ }\bibfield
  {title} {\enquote {\bibinfo {title} {{Continuous Markovian model for L\'evy}
  random walks with superdiffusive and superballistic regimes},}\ }\href
  {\doibase 10.1140/epjb/e2010-10422-4} {\bibfield  {journal} {\bibinfo
  {journal} {The European Physical Journal B}\ }\textbf {\bibinfo {volume}
  {78}},\ \bibinfo {pages} {207--216} (\bibinfo {year} {2010})}\BibitemShut
  {NoStop}%
\bibitem [{Note2()}]{Note2}%
  \BibitemOpen
  \bibinfo {note} {The cause of this issue lies in the fact that $x=0$ is both
  a stable attraction point for the unpertubed motion ($\epsilon \tau ^2=0$)
  and a point where the perturbation vanishes. Hence, given that any perturbed
  trajectory will eventually arrive at $x=0$, this is an accumulation point of
  any initial ``ensemble''.}\BibitemShut {Stop}%
\bibitem [{\citenamefont {Bianucci}\ and\ \citenamefont
  {Grigolini}(1992)}]{bgJCP96}%
  \BibitemOpen
  \bibfield  {author} {\bibinfo {author} {\bibfnamefont {M.}~\bibnamefont
  {Bianucci}}\ and\ \bibinfo {author} {\bibfnamefont {P.}~\bibnamefont
  {Grigolini}},\ }\bibfield  {title} {\enquote {\bibinfo {title} {{Nonlinear
  and non Markovian fluctuation-dissipation processes: A {Fokker-Planck}
  treatment}},}\ }\href {\doibase 10.1063/1.462657} {\bibfield  {journal}
  {\bibinfo  {journal} {The Journal of Chemical Physics}\ }\textbf {\bibinfo
  {volume} {96}},\ \bibinfo {pages} {6138--6148} (\bibinfo {year}
  {1992})}\BibitemShut {NoStop}%
\bibitem [{\citenamefont {Bianucci}\ and\ \citenamefont
  {Mannella}(2020)}]{bmJPC4}%
  \BibitemOpen
  \bibfield  {author} {\bibinfo {author} {\bibfnamefont {M.}~\bibnamefont
  {Bianucci}}\ and\ \bibinfo {author} {\bibfnamefont {R.}~\bibnamefont
  {Mannella}},\ }\bibfield  {title} {\enquote {\bibinfo {title} {Optimal {FPE}
  for non-linear 1{D}-{SDE}. i: Additive {Gaussian} colored noise},}\ }\href
  {\doibase 10.1088/2399-6528/abc54e} {\bibfield  {journal} {\bibinfo
  {journal} {Journal of Physics Communications}\ }\textbf {\bibinfo {volume}
  {4}},\ \bibinfo {pages} {105019} (\bibinfo {year} {2020})}\BibitemShut
  {NoStop}%
\bibitem [{\citenamefont {Bianucci}\ \emph
  {et~al.}(1993{\natexlab{a}})\citenamefont {Bianucci}, \citenamefont {Bonci},
  \citenamefont {Trefan}, \citenamefont {West},\ and\ \citenamefont
  {Grigolini}}]{bbgwgPLA174}%
  \BibitemOpen
  \bibfield  {author} {\bibinfo {author} {\bibfnamefont {M.}~\bibnamefont
  {Bianucci}}, \bibinfo {author} {\bibfnamefont {L.}~\bibnamefont {Bonci}},
  \bibinfo {author} {\bibfnamefont {G.}~\bibnamefont {Trefan}}, \bibinfo
  {author} {\bibfnamefont {B.~J.}\ \bibnamefont {West}}, \ and\ \bibinfo
  {author} {\bibfnamefont {P.}~\bibnamefont {Grigolini}},\ }\bibfield  {title}
  {\enquote {\bibinfo {title} {Brownian motion generated by a two-dimensional
  mapping},}\ }\href {\doibase http://dx.doi.org/10.1016/0375-9601(93)90194-5}
  {\bibfield  {journal} {\bibinfo  {journal} {Physics Letters A}\ }\textbf
  {\bibinfo {volume} {174}},\ \bibinfo {pages} {377 -- 383} (\bibinfo {year}
  {1993}{\natexlab{a}})}\BibitemShut {NoStop}%
\bibitem [{\citenamefont {Bianucci}\ \emph
  {et~al.}(1993{\natexlab{b}})\citenamefont {Bianucci}, \citenamefont
  {Mannella}, \citenamefont {Fan}, \citenamefont {Grigolini},\ and\
  \citenamefont {West}}]{bmxgwPRE47}%
  \BibitemOpen
  \bibfield  {author} {\bibinfo {author} {\bibfnamefont {M.}~\bibnamefont
  {Bianucci}}, \bibinfo {author} {\bibfnamefont {R.}~\bibnamefont {Mannella}},
  \bibinfo {author} {\bibfnamefont {X.}~\bibnamefont {Fan}}, \bibinfo {author}
  {\bibfnamefont {P.}~\bibnamefont {Grigolini}}, \ and\ \bibinfo {author}
  {\bibfnamefont {B.~J.}\ \bibnamefont {West}},\ }\bibfield  {title} {\enquote
  {\bibinfo {title} {Standard fluctuation-dissipation process from a
  deterministic mapping},}\ }\href {\doibase 10.1103/PhysRevE.47.1510}
  {\bibfield  {journal} {\bibinfo  {journal} {Phys. Rev. E}\ }\textbf {\bibinfo
  {volume} {47}},\ \bibinfo {pages} {1510--1519} (\bibinfo {year}
  {1993}{\natexlab{b}})}\BibitemShut {NoStop}%
\bibitem [{\citenamefont {Bianucci}, \citenamefont {West},\ and\ \citenamefont
  {Grigolini}(1994)}]{bwgPLA190}%
  \BibitemOpen
  \bibfield  {author} {\bibinfo {author} {\bibfnamefont {M.}~\bibnamefont
  {Bianucci}}, \bibinfo {author} {\bibfnamefont {B.~J.}\ \bibnamefont {West}},
  \ and\ \bibinfo {author} {\bibfnamefont {P.}~\bibnamefont {Grigolini}},\
  }\bibfield  {title} {\enquote {\bibinfo {title} {Probing microscopic chaotic
  dynamics by observing macroscopic transport processes},}\ }\href {\doibase
  http://dx.doi.org/10.1016/0375-9601(94)90731-5} {\bibfield  {journal}
  {\bibinfo  {journal} {Physics Letters A}\ }\textbf {\bibinfo {volume}
  {190}},\ \bibinfo {pages} {447 -- 454} (\bibinfo {year} {1994})}\BibitemShut
  {NoStop}%
\bibitem [{\citenamefont {Bianucci}(2004)}]{bPRE70}%
  \BibitemOpen
  \bibfield  {author} {\bibinfo {author} {\bibfnamefont {M.}~\bibnamefont
  {Bianucci}},\ }\bibfield  {title} {\enquote {\bibinfo {title} {Ordinary
  chemical reaction process induced by a unidimensional map},}\ }\href
  {\doibase 10.1103/PhysRevE.70.026107} {\bibfield  {journal} {\bibinfo
  {journal} {Phys. Rev. E}\ }\textbf {\bibinfo {volume} {70}},\ \bibinfo
  {pages} {026107--1--12617--6} (\bibinfo {year} {2004})}\BibitemShut {NoStop}%
\bibitem [{\citenamefont {Zwanzig}(2001)}]{Zwanzig2001}%
  \BibitemOpen
  \bibinfo {editor} {\bibfnamefont {R.}~\bibnamefont {Zwanzig}},\ ed.,\ \href
  {http://ukcatalogue.oup.com/product/9780195140187.do} {\emph {\bibinfo
  {title} {Nonequilibrium Statistical Mechanics}}}\ (\bibinfo  {publisher}
  {Oxford University Press},\ \bibinfo {address} {Oxford},\ \bibinfo {year}
  {2001})\BibitemShut {NoStop}%
\bibitem [{\citenamefont {Grigolini}(1989)}]{grigolini1989}%
  \BibitemOpen
  \bibfield  {author} {\bibinfo {author} {\bibfnamefont {P.}~\bibnamefont
  {Grigolini}},\ }\bibfield  {title} {\enquote {\bibinfo {title} {The
  projection approach to the {Fokker-Planck} equation: applications to
  phenomenological stochastic equations with colored noises},}\ }in\ \href
  {\doibase http://dx.doi.org/10.1017/CBO9780511897818} {\emph {\bibinfo
  {booktitle} {Noise in Nonlinear Dynamical Systems}}},\ Vol.~\bibinfo {volume}
  {1},\ \bibinfo {editor} {edited by\ \bibinfo {editor} {\bibfnamefont
  {F.}~\bibnamefont {Moss}}\ and\ \bibinfo {editor} {\bibfnamefont {P.~V.~E.}\
  \bibnamefont {McClintock}}}\ (\bibinfo  {publisher} {Cambridge University
  Press},\ \bibinfo {address} {Cambridge, England},\ \bibinfo {year} {1989})\
  Chap.~\bibinfo {chapter} {5}, p.\ \bibinfo {pages} {161}\BibitemShut
  {NoStop}%
\bibitem [{\citenamefont {Bianucci}\ \emph {et~al.}(1995)\citenamefont
  {Bianucci}, \citenamefont {Mannella}, \citenamefont {West},\ and\
  \citenamefont {Grigolini}}]{bmwgPRE51}%
  \BibitemOpen
  \bibfield  {author} {\bibinfo {author} {\bibfnamefont {M.}~\bibnamefont
  {Bianucci}}, \bibinfo {author} {\bibfnamefont {R.}~\bibnamefont {Mannella}},
  \bibinfo {author} {\bibfnamefont {B.~J.}\ \bibnamefont {West}}, \ and\
  \bibinfo {author} {\bibfnamefont {P.}~\bibnamefont {Grigolini}},\ }\bibfield
  {title} {\enquote {\bibinfo {title} {From dynamics to thermodynamics: Linear
  response and statistical mechanics},}\ }\href {\doibase
  10.1103/PhysRevE.51.3002} {\bibfield  {journal} {\bibinfo  {journal} {Phys.
  Rev. E}\ }\textbf {\bibinfo {volume} {51}},\ \bibinfo {pages} {3002--3022}
  (\bibinfo {year} {1995})}\BibitemShut {NoStop}%
\bibitem [{\citenamefont {Bianucci}(2015)}]{bJSTAT2015}%
  \BibitemOpen
  \bibfield  {author} {\bibinfo {author} {\bibfnamefont {M.}~\bibnamefont
  {Bianucci}},\ }\bibfield  {title} {\enquote {\bibinfo {title} {On the
  correspondence between a large class of dynamical systems and stochastic
  processes described by the generalized {Fokker-Planck} equation with
  state-dependent diffusion and drift coefficients},}\ }\href {\doibase
  10.1088/1742-5468/2015/05/P05016} {\bibfield  {journal} {\bibinfo  {journal}
  {Journal of Statistical Mechanics: Theory and Experiment}\ }\textbf {\bibinfo
  {volume} {2015}},\ \bibinfo {pages} {P05016} (\bibinfo {year}
  {2015})}\BibitemShut {NoStop}%
\bibitem [{\citenamefont {Bianucci}, \citenamefont {Bologna},\ and\
  \citenamefont {Mannella}(2023)}]{bbmSUB}%
  \BibitemOpen
  \bibfield  {author} {\bibinfo {author} {\bibfnamefont {M.}~\bibnamefont
  {Bianucci}}, \bibinfo {author} {\bibfnamefont {M.}~\bibnamefont {Bologna}}, \
  and\ \bibinfo {author} {\bibfnamefont {R.}~\bibnamefont {Mannella}},\
  }\bibfield  {title} {\enquote {\bibinfo {title} {About the optimal {FPE} for
  non-linear 1d-{SDE} with {Gaussian} noise: the pitfall of the perturbative
  approach},}\ }\href@noop {} {\bibfield  {journal} {\bibinfo  {journal} {J.
  Stat. Phys.}\ } (\bibinfo {year} {2023})}\BibitemShut {NoStop}%
\bibitem [{\citenamefont {Bianucci}(2018)}]{bJMP59}%
  \BibitemOpen
  \bibfield  {author} {\bibinfo {author} {\bibfnamefont {M.}~\bibnamefont
  {Bianucci}},\ }\bibfield  {title} {\enquote {\bibinfo {title} {Using some
  results about the {Lie} evolution of differential operators to obtain the
  {Fokker-Planck} equation for {non-Hamiltonian} dynamical systems of
  interest},}\ }\href {\doibase 10.1063/1.5037656} {\bibfield  {journal}
  {\bibinfo  {journal} {Journal of Mathematical Physics}\ }\textbf {\bibinfo
  {volume} {59}},\ \bibinfo {pages} {053303} (\bibinfo {year} {2018})},\
  \Eprint {http://arxiv.org/abs/https://doi.org/10.1063/1.5037656}
  {https://doi.org/10.1063/1.5037656} \BibitemShut {NoStop}%
\bibitem [{Note3()}]{Note3}%
  \BibitemOpen
  \bibinfo {note} {$\protect \tmspace +\thinmuskip {.1667em}_1F_1\left
  (a;b;z\right ):=\DOTSB \sum@ \slimits@ _{k=0}^{\infty } \protect \frac
  {a_k}{b_k k!}$, where $(x)_n=x (x+1) \protect \ldots (x+n-1)=\Gamma
  (x+n)/\Gamma (x) $ is the Pochhammer symbol.}\BibitemShut {Stop}%
\bibitem [{\citenamefont {Bianucci}\ and\ \citenamefont
  {Bologna}(2020)}]{bbJSTAT4}%
  \BibitemOpen
  \bibfield  {author} {\bibinfo {author} {\bibfnamefont {M.}~\bibnamefont
  {Bianucci}}\ and\ \bibinfo {author} {\bibfnamefont {M.}~\bibnamefont
  {Bologna}},\ }\bibfield  {title} {\enquote {\bibinfo {title} {About the
  foundation of the {Kubo} generalized cumulants theory: a revisited and
  corrected approach},}\ }\href {\doibase 10.1088/1742-5468/ab7755} {\bibfield
  {journal} {\bibinfo  {journal} {Journal of Statistical Mechanics: Theory and
  Experiment}\ }\textbf {\bibinfo {volume} {2020}},\ \bibinfo {pages} {043405}
  (\bibinfo {year} {2020})}\BibitemShut {NoStop}%
\end{thebibliography}
%merlin.mbs aipnum4-1.bst 2010-07-25 4.21a (PWD, AO, DPC) hacked
%Control: key (0)
%Control: author (8) initials jnrlst
%Control: editor formatted (1) identically to author
%Control: production of article title (0) allowed
%Control: page (1) range
%Control: year (1) truncated
%Control: production of eprint (0) enabled
 \newcommand{\noop}[1]{}

\end{document}